\documentclass[11pt]{article}
  
\title{On the combinatorial structure of $0/1$-matrices representing nonobtuse simplices}
\date{\today} 
\author{Jan Brandts and Apo Cihangir} 
       
\usepackage{a4wide,authblk,epsfig,amsmath,amssymb,makeidx,graphics,color,tikz} 

\begin{document}                
            
\newtheorem{Th}{Theorem}[section]          
\newtheorem{Le}[Th]{Lemma}       
\newtheorem{Co}[Th]{Corollary}        
\newtheorem{Pro}[Th]{Proposition}            
\newtheorem{Def}[Th]{Definition}          
\newtheorem{rem}[Th]{Remark}              
\newcommand{\be}{\begin{equation}}         
\newcommand{\ee}{\end{equation}}          
\newcommand{\RR}{\mathbb{R}}       
\newcommand{\half}{\frac{1}{2}}  
\newcommand{\hdrie}{\hspace{3mm}}  
\newcommand{\und}{\hdrie\mbox{\rm and }\hdrie}  
\newcommand{\sth}{\hdrie | \hdrie}
    
\def\zeros {{\bf 0}}  
\def\ones{{\bf 1}}
\def\CC {\mathcal{C}}
\def\Bc {\mathcal{B}}
\def\PP {\mathcal{P}}
\def\Bn {\mathcal{B}_n} 
\def\bb {\mathbb{b}}
\def\GG {\mathcal{G}}
\def\OO {\mathcal{O}} 
\def\II {\mathcal{I}} 
\def\BB {\mathbb{B}}  
\def\MM {\mathcal{M}}
\def\EE {\mathcal{E}}
\def\PP {\mathcal{P}}
\def\NN {\mathcal{N}} 
\def\SS {\mathbb{S}}
\def\ZZ {\mathcal{Z}}
\def\TT {\mathcal{T}}
\def\AA {\mathcal{A}}
\def\HH {\mathcal{H}}
\def\csim {\overset{c}{\sim}} 
\def\psim {\overset{p}{\sim}} 
\def\zsim {\overset{{\small 0/1}}{\sim}} 
\def\sort {{\rm sort}}
\def\Eqv {\hdrie\Leftrightarrow\hdrie} 
\def\Bnk {\BB^{n\times k}}
\def\Bnn {\BB^{n\times n}}
\def\Mnn {\MM^{n\times n}}
\def\Mnk {\MM^{n\times k}}  
\def\mod {\,{\rm mod }\,}
\def\conv {\mbox{\rm conv}}
\newcommand{\supp}{{\rm supp}} 
\newcommand{\ol}{\overline} 
\newcommand{\vep}{\varepsilon}
\newcommand{\zok}{$0/1$-$k$} 
\newcommand{\zo}{$0/1$} 
\newcommand{\nmo}{n\!-\!1}
\newcommand{\npo}{n\!+\!1}
\newcommand{\npt}{n\!+\!2}
\def\ptpi {(P^\top P)^{-1}}
\def\Bnk {\BB^{n\times k}}
\def\Bnn {\BB^{n\times n}}
\def\Mnn {\MM^{n\times n}}
\def\Mnk {\MM^{n\times k}} 
\def\conv {\mbox{\rm conv}}
\def\rank {\mbox{\rm rank}}
\maketitle

\begin{abstract}

A {\em $0/1$-simplex} is the convex hull of $n+1$ affinely independent vertices of the 
unit $n$-cube $I^n$. It is {\em nonobtuse} if none its dihedral angles is obtuse, and {\em acute} if additionally 
none of them is right. Acute $0/1$-simplices in $I^n$ can be represented by  $0/1$-matrices $P$ of size 
$n\times n$ whose Gramians $G=P^\top P$ have an inverse that is {\em strictly} diagonally dominant, with {\em negative} off-diagonal entries. 

In this paper, we will prove that the positive part $D$ of the transposed inverse $P^{-\top}$ of 
$P$  is {\em doubly stochastic} and has the same support as $P$. In fact, $P$ has a {\em fully indecomposable} 
doubly stochastic pattern. The negative part $C$ of $P^{-\top}$ is strictly row-substochastic and its support is 
complementary to that of $D$, showing that $P^{-\top}=D-C$ has no zero entries and has positive row sums. 
As a consequence, for each facet $F$ of an acute $0/1$-facet $S$ there exists {\em at most one} other acute 
$0/1$-simplex $\hat{S}$ in $I^n$ having $F$ as a facet. We call $\hat{S}$ the {\em acute neighbor} of $S$ at $F$.

If $P$ represents a $0/1$-simplex that is merely nonobtuse, the inverse of $G=P^\top P$ is only {\em weakly} 
diagonally dominant and has {\em nonpositive} off-diagonal entries. Consequently, $P^{-\top}$ can have 
entries equal to zero. We show that its positive part $D$ is still doubly stochastic, but its support may 
be {\em strictly contained} in the support of $P$. This allows $P$ to have no doubly stochastic pattern 
and to be {\em partly decomposable}. In theory, this might cause a nonobtuse $0/1$-simplex $S$ to 
have several {\em nonobtuse} neighbors $\hat{S}$ at each of its facets.

In the remainder of the paper, we study nonobtuse $0/1$-simplices $S$ having a partly decomposable 
matrix representation $P$. We prove that if $S$ has such a matrix representation, it also has 
a {\em block diagonal} matrix representation with at least two diagonal blocks. Moreover, {\em all} 
matrix representations of $S$ will then be partly decomposable. This proves that the {\em combinatorial} 
property of having a fully indecomposable matrix representation with doubly stochastic pattern is a 
{\em geometrical} property of a subclass of nonobtuse $0/1$-simplices, invariant under all $n$-cube 
symmetries. We will show that a nonobtuse simplex with partly decomposable matrix representation 
can be split in mutually orthogonal simplicial facets whose dimensions add up to $n$, and in which 
each facet has a fully indecomposable matrix representation. Using this insight, we are able to extend 
the {\em one neighbor theorem} for acute simplices to a larger class of nonobtuse simplices.

\end{abstract}

\section{Introduction} 
\begin{figure}[h]
\begin{center}
\begin{tikzpicture}[scale=0.75, every node/.style={scale=0.75}]
\begin{scope}[shift={(0,0)}]
\draw[fill=gray!20!white] (0,0)--(2,0)--(0.5,3)--cycle;
\draw (0,0)--(2,0)--(2,2)--(0,2)--cycle;
\draw (0.5,1)--(2.5,1)--(2.5,3)--(0.5,3)--cycle;
\draw (0,0)--(0.5,1);
\draw (2,0)--(2.5,1);
\draw (0,2)--(0.5,3);
\draw (2,2)--(2.5,3);
\draw (2,0)--(0.5,1);
\draw[fill=black] (0,0) circle [radius=0.05];
\draw[fill=black] (0.5,1) circle [radius=0.05];
\draw[fill=black] (2,0) circle [radius=0.05];
\draw[fill=black] (0.5,3) circle [radius=0.05];
\node[scale=0.8] at (3.4,0.5) {$\left[\begin{array}{rrr} 0 & 0 & 1 \\  1 & 1 & 0\\  1 & 0 & 0\end{array}\right]$};
\end{scope}
\begin{scope}[shift={(4.5,0)}]
\draw[fill=gray!20!white] (0,0)--(2,0)--(2.5,1)--(0,2)--cycle;
\draw (0,0)--(2,0)--(2,2)--(0,2)--cycle;
\draw (0.5,1)--(2.5,1)--(2.5,3)--(0.5,3)--cycle;
\draw (0,0)--(0.5,1);
\draw (2,0)--(2.5,1);
\draw (0,2)--(0.5,3);
\draw (2,2)--(2.5,3);
\draw (2,0)--(0,2);
\draw (0,0)--(2.5,1);
\draw[fill=black] (0,0) circle [radius=0.05];
\draw[fill=black] (2.5,1) circle [radius=0.05];
\draw[fill=black] (2,0) circle [radius=0.05];
\draw[fill=black] (0,2) circle [radius=0.05];
\node[scale=0.8] at (3.4,0.5) {$\left[\begin{array}{rrr} 1 & 1 & 0 \\  0 & 1 & 0\\  0 & 0 & 1\end{array}\right]$};
\end{scope}
\begin{scope}[shift={(9,0)}]
\draw[fill=gray!20!white] (0,0)--(2,0)--(2,2)--(0.5,1)--cycle;
\draw (0,0)--(2,0)--(2,2)--(0,2)--cycle;
\draw (0.5,1)--(2.5,1)--(2.5,3)--(0.5,3)--cycle;
\draw (0,0)--(0.5,1);
\draw (2,0)--(2.5,1);
\draw (0,2)--(0.5,3);
\draw (2,2)--(2.5,3);
\draw (2,0)--(0.5,1);
\draw (0,0)--(2,2);
\draw[fill=black] (0,0) circle [radius=0.05];
\draw[fill=black] (2,2) circle [radius=0.05];
\draw[fill=black] (2,0) circle [radius=0.05];
\draw[fill=black] (0.5,1) circle [radius=0.05];
\node[scale=0.8] at (3.4,0.5) {$\left[\begin{array}{rrr} 1 & 0 & 1 \\ 0 & 1 & 0\\  0 & 0 & 1\end{array}\right]$};
\end{scope}
\begin{scope}[shift={(13.5,0)}]
\draw[fill=gray!20!white] (0,0)--(0,2)--(2.5,3)--(2,2)--cycle;
\draw (0,0)--(2,0)--(2,2)--(0,2)--cycle;
\draw (0.5,1)--(2.5,1)--(2.5,3)--(0.5,3)--cycle;
\draw (0,0)--(0.5,1);
\draw (2,0)--(2.5,1);
\draw (0,2)--(0.5,3);
\draw (2,2)--(2.5,3);
\draw (0,0)--(2,2);
\draw (0,0)--(2.5,3);
\draw[fill=black] (0,0) circle [radius=0.05];
\draw[fill=black] (0,2) circle [radius=0.05];
\draw[fill=black] (2.5,3) circle [radius=0.05];
\draw[fill=black] (2,2) circle [radius=0.05];
\node[scale=0.8] at (3.4,0.5) {$\left[\begin{array}{rrr} 0 & 1 & 1 \\0 & 1 & 0\\ 1 & 1 & 1\end{array}\right]$};
\end{scope}
\end{tikzpicture}
\end{center}
\caption{\small{Matrix representations of the same $0/1$-tetrahedron modulo the action of $\mathcal{B}_3$.}}
\label{Nfigure1}
\end{figure}
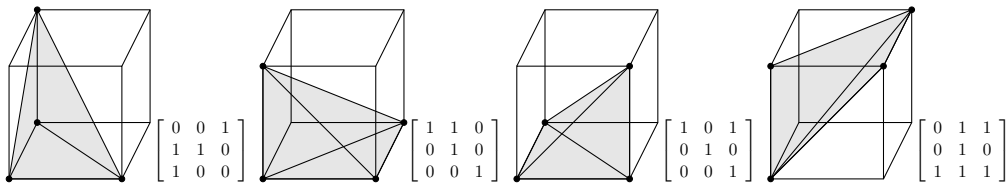

A $0/1$-simplex is an $n$-dimensional $0/1$-polytope \cite{KaZi} with $n+1$ vertices. 
Equivalently, it is the convex hull of $n+1$ of the $2^n$ elements of the set $\BB^n$ of 
vertices of the unit $n$-cube $I^n$ whenever this hull has dimension $n$. Throughout 
this paper, we will study $0/1$-simplices modulo the action of the hyperocthedral group 
$\Bn$ of symmetries of $I^n$. As a consequence, we may assume without loss of 
generality that a $0/1$-simplex $S$ has the origin as a vertex. This makes it possible 
to represent $S$ by a non-singular $n\times n$ matrix $P$ whose columns are the 
remaining $n$ vertices of $S$. Of course, this representation is far from unique, as is 
illustrated by the $0/1$-tetrahedron in Figure~\ref{Nfigure1}. First of all, there is a 
choice which vertex of $S$ is located at the origin. Secondly, column permutations of $P$ 
correspond to relabeling of the nonzero vertices of $S$, and thirdly, row 
permutations correspond to relabeling of the coordinate axis.

\smallskip
  
We will be studying $0/1$-simplices with certain geometric properties. These will be 
invariant under congruence, and in particular invariant under the action of $\Bn$, 
which forms a subset of the congruences of $I^n$. Thus, each of the matrix representations 
carries the required geometric information of the $0/1$-simplex it represents. To be more 
specific, we will study the set of {\em acute} $0/1$-simplices, whose dihedral angles are 
all acute, the {\em nonobtuse} $0/1$-simplices, none of whose dihedral angles is obtuse, 
and the set of {\em orthogonal} simplices. An orthogonal simplex is a nonobtuse simplex 
with exactly $n$ acute dihedral angles and $\half n(n-1)$ right dihedral angles.\\[2mm]
It is not difficult to establish that a $0/1$-simplex $S$ is nonobtuse if and only if the 
inverse $\ptpi$ of the Gramian of any matrix representation $P$ of $S$ is a diagonally 
dominant Stieltjes matrix. This Gramian is {\em strictly} diagonally dominant and has 
even negative off-diagonal entries if and only if $S$ is acute. See \cite{BrCi,BrCi2,BrKoKr} 
for details. In this paper we will study the properties of the $0/1$-matrices 
that represent acute, nonobtuse, and orthogonal simplices.

\subsection{Motivation} 
The motivation to study nonobtuse simplices goes back to their appearance in 
finite element methods \cite{Bra,Bre} to approximate solutions of PDEs, in 
which triangulations consisting of nonobtuse simplices can be used to guarantee 
discrete maximum and comparison principles \cite{BrKoKr2}. We then found that 
they figure in other applications, see \cite{BrKoKrSo} and the references therein. 
In the context of $0/1$-simplices and $0/1$-matrices, it is well known \cite{Gr} 
that the {\em Hadamard conjecture} \cite{Hada} is equivalent to the existence 
of a {\em regular $0/1$-simplex} in each $n$ cube with $n-3$ divisible by $4$. 
Note that a regular simplex is always acute. Thus, studying acute $0/1$-simplices 
can be seen as an attempt to study the Hadamard conjecture in new context, which 
is wider, but not too wide. Indeed, acute $0/1$-simplices, although present in any 
dimension, are still very rare in comparison to {\em all} $0/1$-simplices. 
See \cite{BrCi2}, in which we describe the computational generation of 
acute $0/1$-simplices, as well as several mathematical properties. This 
paper can be seen as a continuation of \cite{BrCi2}, in which some new 
results on acute $0/1$-simplices are presented, as well as on the again 
slightly larger class of nononbtuse $0/1$-simplices.

\subsection{Outline}
We start in Section~\ref{Sect-2} with some preliminaries related to the hyperoctahedral group 
of cube symmetries, to combinatorics, and to the linear algebra behind the geometry of 
nonobtuse and acute simplices. We refer to \cite{BrCi2} for much more detailed information 
on the hyperoctahedral group and combinatorical aspects, and to \cite{BrKoKrSo} for 
applications of nonobtuse simplices. In Section~\ref{Sect-3} we present our new results 
concerning {\em sign properties} of the {\em transposed inverses} $P^{-\top}$ of 
matrix representations $P$ of acute $0/1$-simplices $S$. These results imply that the 
matrices $P$ are {\em fully indecomposable} with {\em doubly stochastic pattern} \cite{Bru1}. 
From this follows the so-called {\em one neighbor theorem}, which states that 
all $(\nmo)$-facets $F$ of $S$ are {\em interior} to the cube, and that each is 
shared by at most one other acute $0/1$-simplex in $I^n$. See \cite{BrDiHaKr} 
for an alternative proof of that fact. If $S$ is merely a nonobtuse $0/1$-simplex, 
the support of $P$ only {\em contains} a doubly stochastic pattern, and moreover, 
$P$ can be {\em partly decomposable}. In Section~\ref{Sect-4} we study the matrix 
representations of such nonobtuse $0/1$-simplices with partly decomposable matrix 
representations. The main conclusion is that each of them consists of $p$ with 
$2\leq p \leq n$ mutually orthogonal facets $F_1,\dots,F_p$ of respective 
dimensions $k_1,\dots,k_p$ that add up to $n$. Moreover, each $k_j\times k_j$ matrix 
representation of each facet $F_j$ is fully indecomposable. In case all facets 
$F_1,\dots,F_n$ are one-dimensional, the corresponding $0/1$-simplex is a 
so-called {\em orthogonal} simplex, as it has a spanning tree of mutually 
orthogonal edges. Orthogonal $0/1$-simplices played an important role in 
the nonobtuse cube triangulation problem, solved in \cite{BrDiHaKr}. 
We briefly recall them in Section~\ref{Sect-5} and put them into the 
novel context of Section~\ref{Sect-4}. Finally, in Section~\ref{Sect-6} 
we use the insights developed so far to prove a one neighbor 
theorem for a wider class of nonobtuse simplices.

\section{Preliminaries}\label{Sect-2}
Let $\BB=\{0,1\}$, and write $\BB^n=\BB^{n\times 1}$ for the 
set of vertices of the unit $n$-cube $I^n=[0,1]^n$, which also 
contains the standard basis vectors $e_1^n,\dots,e_n^n$ and their 
sum $e^n$, the {\em all-ones vector}. The $0/1$-matrices of size 
$n\times k$ we denote by $\BB^{n\times k}$. For each 
$X\in\BB^{n\times k}$, define its {\em antipode} $\ol{X}$ by
\be \ol{X}=e^n(e^k)^\top-X,\ee
and write
\be \ones(X) = (e^n)^\top X e^k \und \zeros(X) = \ones(\ol{X})\ee
for the number of entries of $X$ equal to one, and equal to zero, 
respectively. For any $n\times k$ matrix $X$ define its {\em support} $\supp(X)$ by
\be \supp(X) = \{(i,j)\in\{1,\dots,n\}\times\{1,\dots,k\} \sth (e_i^n)^\top X e_j^k \not=0\}.\ee
We now recall some concepts from combinatorial matrix theory \cite{Bru1,BrRy}.

\begin{Def}{\rm A nonnegative matrix $A$ has a {\em doubly stochastic pattern} 
if there exists a doubly stochastic matrix $D=(d_{ij})$ such that $\supp(D)=\supp(A)$}.
\end{Def}

\begin{Def}{\rm A matrix $A\in\Bnn$ is {\em partly decomposable} if there 
exists a $k\in\{1,\dots,n-1\}$ and permutation matrices $\Pi_1,\Pi_2$ such that 
\be\label{one-1} \Pi_1^\top A \Pi_2 = \left[\begin{array}{cc}A_{11} & A_{12} \\ 0 & A_{22}\end{array}\right], \ee
where $A_{11}$ is a $k\times k$ matrix and $A_{22}$ an $(n-k)\times(n-k)$ 
matrix. If $\Pi_1$ and $\Pi_2$ can be taken equal in (\ref{one-1}) then $A$ 
is called {\em reducible}. If $A$ is not partly decomposable it is called 
{\em fully indecomposable}. If $A$ is not reducible it is called {\em irreducible}.}
\end{Def}  
Note that $A\in\Bnn$ is partly decomposable if and only if there exist 
nonzero $v,w\in\BB^n$ with $\ones(v)+\ones(w)=n$ such that 
$v^\top A w=0$, and that $A$ is reducible if additionally, $w=\ol{v}$.
 
\begin{Le}\label{lem-1} Let $X\in\Bnn$ be nonsingular and $X=[\,X_1\,|\,X_2\,]$ 
a block partition of $X$ where $X_1\in\BB^{n\times k}$ and $X_2\in\BB^{n\times(n-k)}$ 
for some $1\leq k<n$. Suppose that $X_1^\top X_2=0$, and hence,
\be X^\top X = \left[\begin{array}{rl} X_1^\top X_1 & 0 \\ 0 & X_2^\top X_2\end{array}\right]. \ee
Then there exists a permutation $\Pi$ such that
\[ \Pi X = \left[\begin{array}{rl} X_{11} & 0 \\ 0 & X_{22}\end{array}\right], \]
where $X_{11}\in \BB^{k\times k}$ and $X_{22}\in\BB^{(n-k)\times(n-k)}$ are nonsingular.
\end{Le}
{\bf Proof. } If $X_1^\top X_2=0$ then in particular $(e^k)^\top X_1^\top X_2 e^{n-k}=0$, hence
\be \supp(X_1e^k) \cap \supp(X_2e^n) = \emptyset. \ee
Since the rank of $X_1$ equals $k$, the support of $X_1e^k$ 
consists of at least $k$ indices. Similarly, the support of $X_2e^{n-k}$ 
consists of at least $n-k$ indices. Thus, they consist of exactly $k$ and 
$n-k$ indices. Now, let $\Pi$ be a permutation that maps the support of 
$X_1e^k$ onto $\{1,\dots,k\}$, then $\Pi X$ has the required form. \hfill $\Box$
  
\subsection{Matrix representations of $0/1$-simplices modulo $0/1$-equivalence}
The convex hull of any subset $\PP\subset\BB^n$ is called a {\em $0/1$-polytope}. 
As a first simple and intuitive result, we explicitly show that no two 
distinct subsets of $\BB^n$ define the same $0/1$-polytope. 

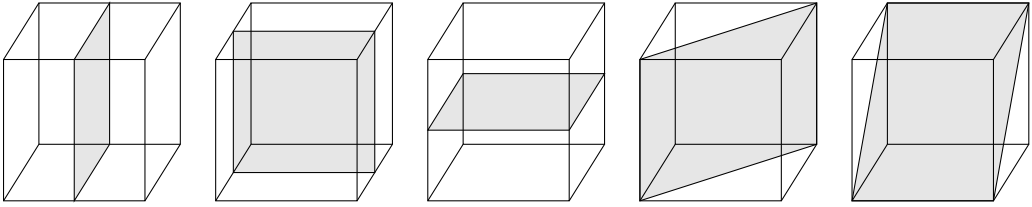
\begin{figure}[h]
\begin{center}
\begin{tikzpicture}[scale=0.94, every node/.style={scale=0.94}]
\begin{scope}
\draw[fill=gray!20!white] (1,0)--(1.5,0.8)--(1.5,2.8)--(1,2)--cycle;
\draw (0,0)--(2,0)--(2,2)--(0,2)--cycle;
\draw (0.5,0.8)--(2.5,0.8)--(2.5,2.8)--(0.5,2.8)--cycle;
\draw (0,0)--(0.5,0.8);
\draw (2,0)--(2.5,0.8);
\draw (0,2)--(0.5,2.8);
\draw (2,2)--(2.5,2.8);
\end{scope}
\begin{scope}[shift={(3,0)}]
\draw[fill=gray!20!white] (0.25,0.4)--(2.25,0.4)--(2.25,2.4)--(0.25,2.4)--cycle;
\draw (0,0)--(2,0)--(2,2)--(0,2)--cycle;
\draw (0.5,0.8)--(2.5,0.8)--(2.5,2.8)--(0.5,2.8)--cycle;
\draw (0,0)--(0.5,0.8);
\draw (2,0)--(2.5,0.8);
\draw (0,2)--(0.5,2.8);
\draw (2,2)--(2.5,2.8); 
\end{scope}
\begin{scope}[shift={(6,0)}]
\draw[fill=gray!20!white] (0,1)--(2,1)--(2.5,1.8)--(0.5,1.8)--cycle;
\draw (0,0)--(2,0)--(2,2)--(0,2)--cycle;
\draw (0.5,0.8)--(2.5,0.8)--(2.5,2.8)--(0.5,2.8)--cycle;
\draw (0,0)--(0.5,0.8);
\draw (2,0)--(2.5,0.8);
\draw (0,2)--(0.5,2.8);
\draw (2,2)--(2.5,2.8); 
\end{scope}
\begin{scope}[shift={(9,0)}]
\draw[fill=gray!20!white] (0,0)--(2.5,0.8)--(2.5,2.8)--(0,2)--cycle;
\draw (0,0)--(2,0)--(2,2)--(0,2)--cycle;
\draw (0.5,0.8)--(2.5,0.8)--(2.5,2.8)--(0.5,2.8)--cycle;
\draw (0,0)--(0.5,0.8);
\draw (2,0)--(2.5,0.8);
\draw (0,2)--(0.5,2.8);
\draw (2,2)--(2.5,2.8); 
\end{scope}
\begin{scope}[shift={(12,0)}]
\draw[fill=gray!20!white] (0,0)--(2,0)--(2.5,2.8)--(0.5,2.8)--cycle;
\draw (0,0)--(2,0)--(2,2)--(0,2)--cycle;
\draw (0.5,0.8)--(2.5,0.8)--(2.5,2.8)--(0.5,2.8)--cycle;
\draw (0,0)--(0.5,0.8);
\draw (2,0)--(2.5,0.8);
\draw (0,2)--(0.5,2.8);
\draw (2,2)--(2.5,2.8); 
\end{scope}
\end{tikzpicture}
\end{center}
\caption{\small{The five reflection planes corresponding to the octahedral group $\mathcal{B}_3$.}}
\label{Nfigure2}
\end{figure}

\begin{Le}\label{lem-2} Suppose that $V\in\BB^{n\times k}$ has distinct 
columns and that $v\in\BB^n$ is not a column of $V$. Then $v$ is not a convex combination of the columns of $V$.
\end{Le} 
{\bf Proof. } Suppose that $v$ is not a column of $V$ and that $Vw=v$ for 
some $w\geq 0$. It suffices to show that $w^\top e^k\not=1$. If $v=0$ then 
$V$ has no zero column, hence $w=0$ is the only nonnegative vector solving 
$Vw=v$. If $v$ has an entry equal to one, say $(e_\ell^n)^\top v=1$, then 
$(e^n_\ell)^\top Vw=1$, hence a non-empty subset $\II\subset\{1,\dots,k\}$ 
of the entries of $w$ sums to one. In order to obtain $w^\top e^k=1$, it is 
necessary that $\supp(w)=\II$. Since $Vw\in\BB^n$, this implies that each 
row of $V$ has either only ones, or only zeros, at its entries at positions in 
$\II$. If $\II$ has more than one element, then $V$ has two or more columns 
that are equal, contradicting its definition.\hfill $\Box$\\[2mm]
This lemma proves the bijective correspondence between the power set of the 
$2^n$ vertices of $I^n$ and the $0/1$-polytopes. Since the cardinality $2^{(2^n)}$ 
of this power set is doubly exponential, the following concept makes sense. Two 
$0/1$-polyoptes will be called {\em $0/1$-equivalent} if they can be transformed 
into one another by the action of an element of the {\em hyperoctahedral group} 
$\Bn$ of affine isometries $\BB^n\rightarrow\BB^n$. The group $\Bn$ has 
$n!2^n$ elements and is generated by:\\[2mm]
$\bullet$ the $n$ reflections in the hyperplanes $2x_i=1$ for $i\in\{1,\dots,n\}$, and\\[2mm]
$\bullet$ the $n-1$ reflections in the hyperplanes $x_i=x_{i+1}$ for $i\in\{1,\dots,n-1\}$.\\[2mm]
See Figure~\ref{Nfigure2} for the five reflection planes of $\mathcal{B}_3$ in $I^3$.

\smallskip

Note that a reflection of the second type exchanges the values of the $i$th 
and  $(i+1)$st coordinate. Hence, products of this type can result in any permutation of the coordinates.\\[3mm]
A $0/1$-{\em simplex} is a $0/1$-polytope in $I^n$ with $n+1$ affinely independent vertices. We will study the set $\SS^n$ of $0/1$-
simplices in $I^n$ modulo the action of $\Bn$. Therefore we can without loss of generality assume that one of its vertices is located at the 
origin, after which $S$ can be conveniently represented by any square matrix whose columns are its remaining $n$ vertices. 
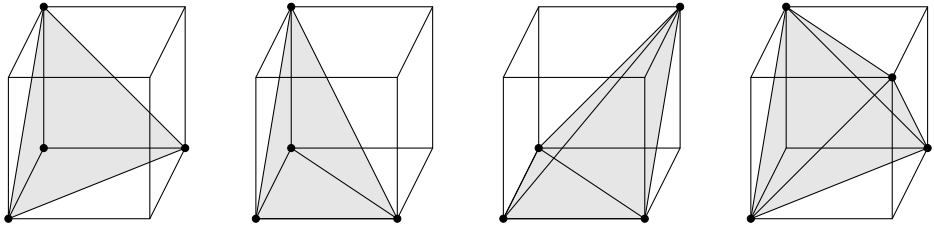
\begin{figure}[h]
\begin{center}
\begin{tikzpicture}[scale=0.94, every node/.style={scale=0.94}]
\begin{scope}
\draw[fill=gray!20!white] (0,0)--(2.5,1)--(0.5,3)--cycle;
\draw (0,0)--(2,0)--(2,2)--(0,2)--cycle;
\draw (0.5,1)--(2.5,1)--(2.5,3)--(0.5,3)--cycle;
\draw (0,0)--(0.5,1);
\draw (2,0)--(2.5,1);
\draw (0,2)--(0.5,3);
\draw (2,2)--(2.5,3);
\draw[fill=black] (0,0) circle [radius=0.05];
\draw[fill=black] (0.5,1) circle [radius=0.05];
\draw[fill=black] (2.5,1) circle [radius=0.05];
\draw[fill=black] (0.5,3) circle [radius=0.05];
\end{scope}
\begin{scope}[shift={(3.5,0)}]
\draw[fill=gray!20!white] (0,0)--(2,0)--(0.5,3)--cycle;
\draw (0,0)--(2,0)--(2,2)--(0,2)--cycle;
\draw (0.5,1)--(2.5,1)--(2.5,3)--(0.5,3)--cycle;
\draw (0,0)--(0.5,1);
\draw (2,0)--(2.5,1);
\draw (0,2)--(0.5,3);
\draw (2,2)--(2.5,3);
\draw (2,0)--(0.5,1);
\draw[fill=black] (0,0) circle [radius=0.05];
\draw[fill=black] (0.5,1) circle [radius=0.05];
\draw[fill=black] (2,0) circle [radius=0.05];
\draw[fill=black] (0.5,3) circle [radius=0.05];
\end{scope}
\begin{scope}[shift={(7,0)}]
\draw[fill=gray!20!white] (0,0)--(0.5,1)--(2.5,3)--(2,0)--cycle;
\draw (0,0)--(2,0)--(2,2)--(0,2)--cycle;
\draw (0.5,1)--(2.5,1)--(2.5,3)--(0.5,3)--cycle;
\draw (0,0)--(0.5,1);
\draw (2,0)--(2.5,1);
\draw (0,2)--(0.5,3);
\draw (2,2)--(2.5,3);
\draw (0,0)--(2.5,3);
\draw (2,0)--(0.5,1);
\draw[fill=black] (0,0) circle [radius=0.05];
\draw[fill=black] (0.5,1) circle [radius=0.05];
\draw[fill=black] (2,0) circle [radius=0.05];
\draw[fill=black] (2.5,3) circle [radius=0.05];
\end{scope}
\begin{scope}[shift={(10.5,0)}]
\draw[fill=gray!20!white] (0,0)--(2.5,1)--(2,2)--(0.5,3)--cycle;
\draw (0,0)--(2,0)--(2,2)--(0,2)--cycle;
\draw (0.5,1)--(2.5,1)--(2.5,3)--(0.5,3)--cycle;
\draw (0,0)--(0.5,1);
\draw (2,0)--(2.5,1);
\draw (0,2)--(0.5,3);
\draw (2,2)--(2.5,3);
\draw (0,0)--(2,2);
\draw (2.5,1)--(0.5,3);
\draw[fill=black] (0,0) circle [radius=0.05];
\draw[fill=black] (0.5,3) circle [radius=0.05];
\draw[fill=black] (2.5,1) circle [radius=0.05];
\draw[fill=black] (2,2) circle [radius=0.05];
\end{scope}
\end{tikzpicture}
\end{center}
\caption{\small{The only four $0/1$-tetrahedra in $I^3$ modulo the action of $\mathcal{B}_3$. In dimension $5$ and higher there exist 
congruent simplices that are not equal modulo $\Bn$.}}
\label{Nfigure3}
\end{figure}  
Modulo the action of $\Bn$, two matrices describe the same $0/1$-simplex if and only if one can be obtained from the other using any 
combination of the following three operations,\\[2mm]
(C) Column permutations;\\[2mm]
(R) Row permutations;\\[2mm]
(X) Select a column $c$ and replace each remaining column $d$ by $(c+d)(\mod 2)$.\\[2mm]
Operation (C) does not correspond to any action of $\Bn$, it just relabels the vertices of $S$. Operations of type (R) correspond to products of 
reflections in hyperplanes $x_i=x_{i+1}$. The operation (X) corresponds to reflecting the vertex $c$ to the origin (and the origin to $c$). It is 
indicated by the letter X because the corresponding matrix operation can also be interpreted as taking the logical {\em exclusive-or} operation 
between a fixed column and the remaining ones.

\subsection{Dihedral angles of $0/1$-simplices}
Given a vertex $v$ of a simplex $S\in\SS^n$, the convex hull of the remaining vertices of $S$ is the {\em facet} $F_v$ of $S$ opposite $v$. 
The {\em dihedral} angle $\alpha$ between two given facets of $S$ is the angle supplementary to the angle $\gamma$ between two normal 
vectors to those facets, both pointing into $S$ or both pointing out of $S$. In other words, $\alpha+\gamma=\pi$.

\begin{Def}{\rm  A simplex $S\in\SS^n$ will be called {\em nonobtuse} if none of its dihedral angles is obtuse (greater than $\pi/2$), and 
{\em acute} if all its dihedral angles are acute (less than $\pi/2$).}
\end{Def} 
If $P$ is a matrix representation of $S$ with columns $p_1,\dots,p_n$, then the columns $q_1,\dots,q_n$ of the matrix $Q=P^{-\top}$ are 
inward normals to the facets $F_{p_1},\dots,F_{p_n}$, respectively, as $Q^\top P=I$.
The vector $q$ satisfying $P^\top q=e^n$ is orthogonal to each difference of two columns of $P$. It is an outward normal to the facet 
$F_{p_0}$ opposite the origin $p_0$ and equals $q=P^{-\top}e^n = q_1+\dots+q_n$. This proves the following proposition.

\begin{Pro}[\cite{BrCi,BrKoKr}]\label{pro-1} Let $F_1,\dots,F_n$ be the facets of $S\in\SS^n$ meeting at the origin and $F_0$ its facet 
opposite the origin, and let $P$ be a matrix representation of $S$. Then $F_0$ makes nonobtuse dihedral angles with $F_1,\dots,F_n$ if and 
only if for all $i\in\{1,\dots,n\}$,
\be\label{eq-1.1}  (e_i^n)^\top(P^\top P)^{-1}e^n \geq 0. \ee
Moreover, each pair of facets $F_i$ and $F_j$ with $i\not=0\not=j$ makes a nonobtuse dihedral angle if and only if
\be\label{eq-1.2} (e_i^n)^\top(P^\top P)^{-1}e_j^n \leq 0,  \ee
for all $i,j\in\{1,\dots,n\}$. Therefore, $S$ is nonobtuse if and only if both (\ref{eq-1.1}) and (\ref{eq-1.2}) hold.
\end{Pro}
Property (\ref{eq-1.1}) translates as {\em diagonal dominance} of $(P^\top P)^{-1}$, and (\ref{eq-1.2}) is called the {\em Stieltjes property} 
of $(P^\top P)^{-1}$, which is then called a {\em Stieltjes matrix}.

\begin{figure}[h]
\begin{center}
\begin{tikzpicture}[scale=0.94, every node/.style={scale=0.94}]
\begin{scope}[shift={(0,0)}]
\draw[gray!10!white, fill=gray!10!white] (2.5,0)--(0.5,1)--(0,2.5)--(3,3.5)--cycle;
\draw[gray!30!white, fill=gray!30!white] (0.5,1)--(0,2.5)--(3,3.5)--cycle;
\draw[->] (2.5,0)--(0.5,1);
\draw[->] (2.5,0)--(0,2.5);
\draw[->] (2.5,0)--(3,3.5);
\draw[fill] (2.5,0) circle [radius=0.07];
\node at (2.8,0.1) {$v_0$};
\node at (1.3,2.5) {$F_0$};
\end{scope}
\begin{scope}[shift={(3.7,0)}]
\draw[gray!10!white, fill=gray!10!white] (2.5,0)--(0.5,1)--(0,2.5)--(3,3.5)--cycle;
\draw[gray!30!white, fill=gray!30!white] (2.5,0)--(0,2.5)--(3,3.5)--cycle;
\draw[->] (0.5,1)--(2.5,0);
\draw[->] (0.5,1)--(0,2.5);
\draw[->] (0.5,1)--(3,3.5);
\draw[fill] (0.5,1) circle [radius=0.07];
\node at (0.3,0.7) {$v_1$};
\node at (1.8,1.5) {$F_1$};
\end{scope}
\begin{scope}[shift={(7.4,0)}]
\draw[gray!10!white, fill=gray!10!white] (2.5,0)--(0.5,1)--(0,2.5)--(3,3.5)--cycle;\
\draw[gray!30!white, fill=gray!30!white] (0.5,1)--(2.5,0)--(3,3.5)--cycle;
\draw[->] (0,2.5)--(2.5,0);
\draw[->] (0,2.5)--(0.5,1);
\draw[->] (0,2.5)--(3,3.5);
\draw[fill] (0,2.5) circle [radius=0.07];
\node at (-0.3,2.2) {$v_2$};
\node at (2.1,1.2) {$F_2$};
\end{scope}
\begin{scope}[shift={(11.1,0)}]
\draw[gray!10!white, fill=gray!10!white] (2.5,0)--(0.5,1)--(0,2.5)--(3,3.5)--cycle;
\draw[gray!30!white, fill=gray!30!white] (0.5,1)--(0,2.5)--(2.5,0)--cycle;
\draw[->] (3,3.5)--(2.5,0);
\draw[->] (3,3.5)--(0.5,1);
\draw[->] (3,3.5)--(0,2.5);
\draw[fill] (3,3.5) circle [radius=0.07];
\node at (3.3,3.2) {$v_3$};
\node at (1,1) {$F_3$};
\end{scope}
\end{tikzpicture}
\end{center}
\caption{\small{Dihedral angles are present in different ways in different matrix representations.}}
\label{Nfigure4}
\end{figure}
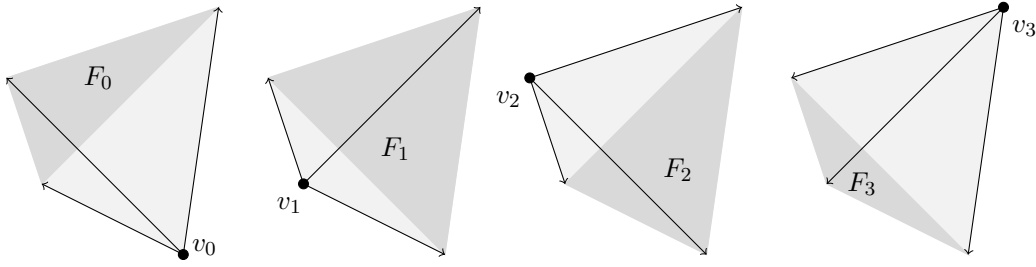

\begin{rem}\label{rem-1}{\rm Condition (\ref{eq-1.1}) is equivalent with the statement that the vertex $v$ of $S$ at the origin orthogonally 
projects onto its opposite facet $F_v$. This proves that the following four statements are equivalent:\\[2mm]
$\bullet$ $S$ is nonobtuse;\\[2mm]
$\bullet$ each vertex $v$ of $S$ projects orthogonally onto its opposite facet $F_v$;\\[2mm]
$\bullet$ each matrix representation $P$ of $S$ satisfies $(P^\top P)^{-1}e^n \geq 0$;\\[2mm]
$\bullet$ each matrix representation $P$ of $S$ satisfies $(e_i^n)^\top(P^\top P)^{-1}e_j^n \leq 0$ for all $i,j\in\{1,\dots,n\}$.}
\end{rem} 
In fact, the second and third condition in Remark \ref{rem-1} can be slightly relaxed. See Figure~\ref{Nfigure4}.
 
\smallskip    

For each $j\in\{0,1,2,3\}$, let $P_j$ be a matrix representations of a tetrahedron $S\in\SS^3$ with its vertex $v_j$ located at the origin. If for 
instance $(P_j^\top P_j)^{-1}e^n \geq 0$ for $j\in\{1,2,3\}$ then $F_1,F_2$ and $F_3$ make only nonobtuse dihedral angles. These include 
{\em all} the six dihedral angles of $S$.

\begin{rem}\label{rem-2} {\rm To characterize acute simplices similarly, replace $\geq$ in (\ref{eq-1.1}) by $>$ and $\leq$ in (\ref{eq-1.2}) by 
$<$. Moreover, replace {\em onto} by {\em into the interior} of its opposite facet.}
\end{rem}
The following two simple combinatorial lemmas will be used further on in this paper.
 
\begin{Le}\label{lem-3} Let $P\in\Bnn$ represent a nonobtuse $0/1$-simplex $S$. Then for all $v\in\BB^n$, 
\be v^\top(P^\top P)^{-1}\ol{v} \leq 0.\ee
If $S$ is even acute then 
\be v^\top(P^\top P)^{-1}\ol{v} < 0 \ee
for all $v\in\BB^n$ with $v\not\in\{0,e^n\}$.
\end{Le} 
{\bf Proof. } In fact, $v^\top(P^\top P)^{-1}\ol{v}$ is the sum of $\ones(v)\ones(\ol{v})$ of the off-diagonal entries of $(P^\top P)^{-1}$. 
According to Proposition \ref{pro-1} these are nonpositive if $S$ is nonobtuse. According to Remark \ref{rem-2} they are negative if $S$ is 
acute. \hfill $\Box$ 
 
\begin{Le}\label{lem-4} Let $P\in\Bnn$ represent a nonobtuse $0/1$-simplex $S$. If for some $v\in\BB^n$,
\be v^\top (P^\top P)^{-1} \ol{v} = 0, \ee
then $v=0$ or $\ol{v}=0$ or $P^\top P$ is reducible.
\end{Le}
{\bf Proof. } Let $v\not=0\not=\ol{v}$ and write $k=\zeros(v)$. Then $1\leq k \leq n-1$. Let $\Pi$ be a permutation such that 
$\supp(\Pi\ol{v})=\{1,\dots,k\}$. Writing $w=\Pi v$, we have that $\ol{w}=\Pi\ol{v}$ and
\be 0 = v^\top (P^\top P)^{-1} \ol{v} = v^\top\Pi^\top\Pi (P^\top P)^{-1} \Pi^\top\Pi \ol{v} = w\Pi (P^\top P)^{-1} \Pi^\top\ol{w},\ee
which is the sum of the entries of $\Pi (P^\top P)^{-1}\Pi^\top$ with indices $(i,j)$ with $k+1\leq i\leq n$ and $1\leq j \leq k$. Since these 
entries are non-positive and their sum equals zero, they are all zero, leading to an $(n-k)\times k$ bottom left block of zeros in 
$\Pi (P^\top P)^{-1}\Pi^\top$. Thus, $(P^\top P)^{-1}$ is reducible, and hence, so is its inverse $P^\top P$. \hfill $\Box$\\[3mm]
A final important observation is the following classical result by Fiedler.

\begin{Le}[\cite{Fie}]\label{lem-5} All $k$-facets of a nonobtuse simplex are nonobtuse and all $k$-facets of an acute simplex are acute.
\end{Le}
It is well known that the converse does not hold. Simplices whose facets are all nonobtuse or acute were studied recently in \cite{BrCi}. 
 
\section{Doubly stochastic patterns and full indecomposability}\label{Sect-3}
We start our investigations with some results on matrix representations of acute $0/1$-simplices. The first one gives a remarkable connection 
with doubly stochastic matrices. It follows from the observation in Remark \ref{rem-1} that for an acute $0/1$-simplex, the altitude from each 
vertex points into the interior of $I^n$. This, in turn, defines the signs of the entries of the normals to its facets in terms of the supports of 
their corresponding vertices.

\begin{Th}\label{th-1} Let $P\in\Bnn$ be a matrix representation of an acute $0/1$-simplex $S\in\SS^n$, and write $Q=P^{-\top}$. Then 
\be\label{eq-3.1} q_{ij}>0 \Leftrightarrow p_{ij} = 1 \und q_{ij}<0 \Leftrightarrow p_{ij}=0. \ee
Defining $0\leq C=(c_{ij})$ and $0\leq D=(d_{ij})$ by 
\be\label{eq-3.2} C = \half\left( |Q|-Q\right) \und D = \half\left(|Q|+Q\right), \ee
where $|Q|$ is the matrix whose entries are the moduli of the entries of $Q$, we have that 
\be Q = D-C, \ee
where $D$ is doubly stochastic and $C$ row-substochastic. 
\end{Th}  
{\bf Proof}. The $j$th column $q_j$ of $Q$ is an inward normal to the facet $F_j$ of $S$ opposite the $j$th column $p_j$ of $P$. Thus, 
$p_j-\alpha q_j$ is an element of the interior of $I^n$ for $\alpha>0$ small enough. From this, (\ref{eq-3.1}) immediately follows. Combining 
this with the fact that the inner product between $p_j$ and $q_j$ equals one, the positive elements in each column of $Q$ add to one. But 
since $Q^\top P=I$, also the inner products between corresponding rows of $P$ and $Q$ equals one, and thus also the positive elements in 
each row of $Q$ add to one. Finally, $Qe^n$ is the outward normal to the facet of $S$ opposite the origin and thus it points into the interior of 
$I^n$. Consequently, $Qe^n>0$, hence $De^n>Ce^n \geq 0$, which shows that $C$ is row-substochastic. \hfill $\Box$ 
 
\begin{rem}\label{rem-3}{\em For $n\geq 7$ there exist matrix representations $P\in\BB^{n\times n}$ of acute $0/1$-simplices in $I^n$  for 
which the matrix $C$ in (\ref{eq-3.2}) is not column-substochastic. This shows that Theorem \ref{th-1} cannot be strengthened in this 
direction. It also proves that if $P$ represents an acute $0/1$-simplex, its transpose $P^\top$ may not do the same. See Figure~\ref{Nfigure5} for an 
example.}\end{rem}
  
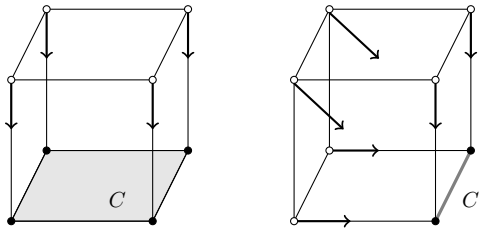
\begin{figure}[h]
\begin{center}
\begin{tikzpicture}[scale=0.94, every node/.style={scale=0.94}]
\draw[fill=gray!20!white] (0,0)--(2,0)--(2.5,1)--(0.5,1)--cycle;
\draw (0,0)--(2,0)--(2,2)--(0,2)--cycle;
\draw (0.5,1)--(2.5,1)--(2.5,3)--(0.5,3)--cycle;
\draw (0,0)--(0.5,1);
\draw (2,0)--(2.5,1);
\draw (0,2)--(0.5,3);
\draw (2,2)--(2.5,3);
\node[scale=0.8] at (1.5,0.3) {$C$};
\draw[fill=black] (0,0) circle [radius=0.05];
\draw[fill=black] (2,0) circle [radius=0.05];
\draw[fill=black] (2.5,1) circle [radius=0.05];
\draw[fill=black] (0.5,1) circle [radius=0.05];
\draw[fill=white] (0,2) circle [radius=0.05];
\draw[fill=white] (2,2) circle [radius=0.05];
\draw[fill=white] (0.5,3) circle [radius=0.05];
\draw[fill=white] (2.5,3) circle [radius=0.05];
\draw[thick,->] (0,1.95)--(0,1.3);
\draw[thick,->] (2,1.95)--(2,1.3);
\draw[thick,->] (0.5,2.95)--(0.5,2.3);
\draw[thick,->] (2.5,2.95)--(2.5,2.3);
\begin{scope}[shift={(4,0)}]
\draw (0,0)--(2,0)--(2,2)--(0,2)--cycle;
\draw (0.5,1)--(2.5,1)--(2.5,3)--(0.5,3)--cycle;
\draw (0,0)--(0.5,1);
\draw[gray,very thick]  (2,0)--(2.5,1);
\draw (0,2)--(0.5,3);
\draw (2,2)--(2.5,3);
\draw[fill=white] (0,0) circle [radius=0.05];
\draw[fill=black] (2,0) circle [radius=0.05];
\draw[fill=black] (2.5,1) circle [radius=0.05];
\draw[fill=white] (0.5,1) circle [radius=0.05];
\draw[fill=white] (0,2) circle [radius=0.05];
\draw[fill=white] (2,2) circle [radius=0.05];
\draw[fill=white] (0.5,3) circle [radius=0.05];
\draw[fill=white] (2.5,3) circle [radius=0.05];
\draw[thick,->] (0,1.95)--(0.7,1.3);
\draw[thick,->] (2,1.95)--(2,1.3);
\draw[thick,->] (0.5,2.95)--(1.2,2.3);
\draw[thick,->] (2.5,2.95)--(2.5,2.3);
\draw[thick,->] (0.05,0)--(0.8,0);
\draw[thick,->] (0.55,1)--(1.2,1);
\node[scale=0.8] at (2.5,0.3) {$C$};
\end{scope}
\end{tikzpicture}
\end{center}
\caption{\small{For each facet $C$ of $I^n$, each vertex of $I^n$ projects onto 
a vertex of $C$. Consequently, no acute $0/1$-simplex has a 
facet contained in a facet of $I^n$.}}
\label{Nfigure6}
\end{figure}

\begin{Co}\label{co-1} Let $P\in\Bnn$ be a matrix representation of an acute $0/1$-simplex $S\in\SS^n$. Then $P$ has a fully 
indecomposable doubly stochastic pattern. 
\end{Co}
{\bf Proof. } Due to (\ref{eq-3.1}), the matrix $D$ in (\ref{eq-3.2}) has the same support as $P$, hence $P$ has a doubly stochastic pattern. 
Next, assume to the contrary that $P$ is partly decomposable, then there exist permutations $\Pi_1,\Pi_2$ such that
\[ \Pi_1^\top P \Pi_2 = \left[\begin{array}{cc}P_{11} & P_{12} \\ 0 & P_{22}\end{array}\right], \]
where $P_{11}$ is a $k\times k$ matrix and $P_{22}$ an $(n-k)\times(n-k)$ matrix for some $k\in\{1,\dots,n-1\}$. But then $Q=P^{-\top}$ 
has entries equal to zero, which contradicts (\ref{eq-3.1}) in Theorem \ref{th-1}. \hfill $\Box$\\[3mm]  
Theorem \ref{th-1}, Corollary \ref{co-1}, and Remark \ref{rem-3} are all illustrated in Figure~\ref{Nfigure5}.
\begin{figure}[h]
\begin{center}
\begin{tikzpicture}[scale=0.94, every node/.style={scale=0.94}]
\draw[fill=gray!20!white] (0,0)--(2.5,1)--(2,2)--cycle;
\draw[gray] (0,0)--(0.5,3)--(2,2);
\draw[thick,->] (1.5,1.2)--(0.6,2.75);
\node[scale=0.9] at (0.2,3) {$p$};
\node[scale=0.9] at (1,1.5) {$q$};
\node[scale=0.8] at (1.2,0.7) {$F_p$};
\node[scale=0.9] at (1.7,1.2) {$\pi$};
\draw[fill=white] (1.5,1.2) circle [radius=0.05];
\draw (0,0)--(2,0)--(2,2)--(0,2)--cycle;
\draw (0.5,1)--(2.5,1)--(2.5,3)--(0.5,3)--cycle;
\draw (0,0)--(0.5,1);
\draw (2,0)--(2.5,1);
\draw (0,2)--(0.5,3);
\draw (2,2)--(2.5,3);
\draw (0,0)--(2,2);
\draw[gray] (2.5,1)--(0.5,3);
\draw[fill=black] (0,0) circle [radius=0.05];
\draw[fill=black] (0.5,3) circle [radius=0.05]; 
\draw[fill=black] (2.5,1) circle [radius=0.05];
\draw[fill=black] (2,2) circle [radius=0.05];
\node[scale=0.73] at (8.5,1.5)  {$P=\left[\begin{array}{rrrrrrr}
     \fbox{1}  &   \fbox{1}  &   \fbox{1}  &   0  &   0  &   \fbox{1}  &   \fbox{1}\\
     \fbox{1}  &   0  &   0  &   \fbox{1}  &   \fbox{1}  &   0  &   0\\
     0  &   \fbox{1}  &   0  &   \fbox{1}  &   \fbox{1}  &   0  &   0\\
     0  &   0  &   \fbox{1}  &   \fbox{1}  &   \fbox{1}  &   \fbox{1}  &   0\\
     0  &   0  &   \fbox{1}  &   \fbox{1}  &   \fbox{1}  &   0  &   \fbox{1}\\
     0  &   0  &   0  &   \fbox{1}  &   0  &   \fbox{1}  &   \fbox{1}\\
     0  &   0  &   0  &   0  &   \fbox{1}  &   \fbox{1}  &   \fbox{1}
     \end{array}\right], \hspace{3mm} P^{-\top}=\dfrac{1}{13}\left[\begin{array}{ccccccc}     
     \fbox{4}  &   \fbox{4}   &  \fbox{3}  &  -2  &  -2  &   \fbox{1}  &   \fbox{1}\\
     \fbox{9}  &  -4   & -3  &   \fbox{2}  &   \fbox{2}  &  -1  &  -1\\
    -4  &   \fbox{9}   & -3  &   \fbox{2}  &   \fbox{2}  &  -1  &  -1\\
    -2  &  -2   &  \fbox{5}  &   \fbox{1}  &   \fbox{1}  &   \fbox{6}  &  -7\\
    -2  &  -2   &  \fbox{5}  &   \fbox{1}  &   \fbox{1}  &  -7  &   \fbox{6}\\
    -1  &  -1   & -4  &   \fbox{7}  &  -6  &   \fbox{3}  &   \fbox{3}\\
    -1  &  -1   & -4  &  -6  &   \fbox{7}  &   \fbox{3}  &   \fbox{3}
    \end{array}\right]$};
\end{tikzpicture}
\end{center}
\caption{\small{In an acute $0/1$-simplex, each vertex $p$ following the inward normal $q$ to its opposite facet $F_p$ (in converse direction) 
projects as $\pi$ in the interior of $F_p$ and hence in the interior of $I^n$. This fixes the signs of the entries of $q$ in terms of those of $p$. 
The matrices $P$ and $P^{-\top}$ (not related to the depicted tetrahedron) constitute an example of the linear algebraic consequences. The 
positions of the positive entries (boxed) of $P^{-\top}$ and $P$ coincide. The positive part $D$ of $P^{-\top}$ is doubly-stochastic. The 
negated negative part $C$ of $P^{-\top}$ is a row-substochastic. It is not column-substochastic because the third column of $C$ adds to 
$\frac{14}{13}$. Thus,  even though also $P^\top$ has a doubly stochastic pattern and is fully indecomposable, it is not a matrix 
representation of an acute binary simplex.}}
\label{Nfigure5}
\end{figure}
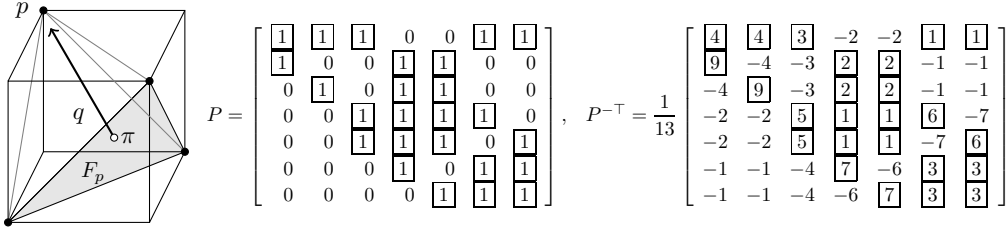  
\begin{rem}{\rm Because {\em each} matrix representation $P$ of an acute $0/1$-simplex has a fully indecomposable doubly stochastic 
pattern, applying to such a matrix $P$ any operation of type (X), as described below Figure~\ref{Nfigure3}, results in another matrix with a fully 
indecomposable doubly stochastic pattern. From the linear algebraic point of view, this is remarkable because generally, both the $0/1$-matrix 
properties of full indecomposability and of having a doubly stochastic pattern are destroyed under operations of type (X). See for instance
\be \small \left[\begin{array}{rrr} 1 & 1 & 1 \\ 1 & 1 & 1 \\ 1 & 1 & 1\end{array}\right] \overset{\mbox{\rm(X)}}{\longrightarrow} \left[\begin{array}{rrr} 1 & 0 & 0 \\ 1 & 0 & 0 \\ 1 & 0 & 0\end{array}\right] \und \left[\begin{array}{rrr} 1 & 0 & 0 \\ 0 & 1 & 0 \\ 0 & 0 & 1\end{array}\right] \overset{\mbox{\rm(X)}}{\longrightarrow} \left[\begin{array}{rrr} 1 & 1 & 1 \\ 0 & 1 & 0 \\ 0 & 0 & 1\end{array}\right]. \ee
From the geometric point of view, this is easy to understand, as the fact that each altitude from each vertex of $S$ points into the interior of 
$I^n$ is invariant under the action of $\Bn$.}
\end{rem} 
The geometric translation of Corollary \ref{co-1} is that if $S\in\SS^n$ is acute, none of its $k$-dimensional facets is contained in a $k$-
dimensional facet of $I^n$ for $k\in\{1,\dots,n-1\}$. The geometric {\em proof} of this is to note that, given a $k$-facet $C$ of $I^n$, no 
vertex of $I^n$ orthogonally projects into {\em the interior} of $C$. In fact, each vertex of $I^n$ projects on a {\em vertex} of $C$. See 
Figure~\ref{Nfigure6}. Thus, if an arbitrary $0/1$-simplex $S$ has a $k$-facet $K$ contained in $C$, each remaining vertex of $S$ projects on a vertex of 
$C$. Remark \ref{rem-2} now shows that $S$ cannot be acute.

\smallskip

Contrary to an acute $0/1$-simplex, a {\em nonobtuse} $0/1$-simplex $S$ may indeed have a $k$-facet $K$ that is contained in a cube facet 
$C$ of $I^n$. If this is the case, then each remaining vertex of $S$ projects on a vertex of $K$. Moreover, $S$ has a partly decomposable 
matrix representation. Before discussing this structure, we first formulate the equivalent of Theorem \ref{th-1} for nonobtuse simplices and 
discuss some of the differences with Theorem \ref{th-1} using an example.

\begin{Th}\label{th-2} Let $P\in\Bnn$ be a matrix representation of a nonobtuse $0/1$-simplex $S\in\SS^n$, and write $Q=P^{-\top}$. Then 
\be\label{eq-3.1-n} q_{ij}\geq 0 \Leftrightarrow p_{ij} = 1 \und q_{ij}\leq 0 \Leftrightarrow p_{ij}=0. \ee
Defining $0\leq C=(c_{ij})$ and $0\leq D=(d_{ij})$ by 
\be\label{eq-3.2-n} C = \half\left( |Q|-Q\right) \und D = \half\left(|Q|+Q\right), \ee
where $|Q|$ is the matrix whose entries are the moduli of the entries of $Q$, we have that 
\be Q = D-C, \ee
where $D$ is doubly stochastic and $C$ row-substochastic. 
\end{Th}  
{\bf Proof}. The proof only differs from the proof of Theorem \ref{th-1} in the sense that $p_j-\alpha q_j$ is now an element of $I^n$ 
including its boundary. This accounts for the $\geq$ and $\leq$ signs in (\ref{eq-3.1-n}) in comparison to the $>$ and $<$ signs in 
(\ref{eq-3.1}). \hfill $\Box$\\[3mm]
Theorem \ref{th-2} is rather weaker than Theorem \ref{th-1}. First of all, the matrix $P^{-\top}$ can have entries equal to zero. Moreover, it 
cannot anymore be concluded that $P$ has a doubly stochastic pattern, only that it {\em contains} a doubly stochastic pattern,
\be \supp(D)\subset \supp(P). \ee
A typical example of this is the following matrix representation $P$ of a nonobtuse $0/1$-simplex,
\begin{figure}[h]
\begin{center}
\begin{tikzpicture}[scale=0.94, every node/.style={scale=0.94}]
\draw[fill=gray!20!white] (0,0)--(2.5,1)--(2.5,3)--cycle;
\draw[gray] (2,0)--(2.5,3);
\draw[thick,->] (1.25,0.5)--(1.9,0.1);
\node[scale=0.9] at (2.2,0) {$p$};
\node[scale=0.9] at (1.3,0.2) {$q$};
\node[scale=0.8] at (1.6,1.3) {$F_p$};
\node[scale=0.9] at (1.4,0.7) {$\pi$};
\draw[fill=white] (1.25,0.5) circle [radius=0.05];
\draw (0,0)--(2,0)--(2,2)--(0,2)--cycle;
\draw (0.5,1)--(2.5,1)--(2.5,3)--(0.5,3)--cycle;
\draw (0,0)--(0.5,1);
\draw (2,0)--(2.5,1);
\draw (0,2)--(0.5,3);
\draw (2,2)--(2.5,3);
\draw[fill=black] (0,0) circle [radius=0.05];
\draw[fill=black] (2,0) circle [radius=0.05]; 
\draw[fill=black] (2.5,1) circle [radius=0.05];
\draw[fill=black] (2.5,3) circle [radius=0.05];
\node[scale=0.73] at (8.5,1.5)  {$P=\left[\begin{array}{ccccccc}
     \fbox{1}  &   \fbox{1}   &    0  &   0  &    \fbox{1}   &    \fbox{1}   &    \fbox{1} \\
     0  &    \fbox{1}   &   0  &   0   &    \fbox{1}   &    \fbox{1}   &    \fbox{1}\\
     0  &   0  &    \fbox{1}   &   \fbox{1}   &    0   &   0   &    0 \\
     0  &   0  &   0  &    \fbox{1}   &   0  &   0  &   0\\
     0  &   0  &   0  &   0  &    \fbox{1}   &    \fbox{1}   &   0\\ 
     0  &   0  &   0  &   0  &    \fbox{1}   &   0  &    \fbox{1} \\
     0  &   0  &   0  &   0  &   0  &    \fbox{1}   &    \fbox{1} 
     \end{array}\right], \hspace{3mm} P^{-\top} = \dfrac{1}{2}
     \left[\begin{array}{rrrrrrr}
      \fbox{2}   &   0  &   0  &   0  &   0  &   0  &   0\\
     -2  &    \fbox{2}   &   0  &   0  &   0  &   0  &   0\\
     0 &   0  &    \fbox{2}   &   0  &   0  &   0  &   0\\
     0  &  0  &   -2  &    \fbox{2}   &   0  &   0  &   0\\
     0  &   -1 &  0  &   0  &    \fbox{1}   &    \fbox{1}   &  -1\\
     0  &   -1 &  0  &   0  &    \fbox{1}   &  -1  &    \fbox{1} \\
     0  &   -1 &  0  &   0  &  -1  &    \fbox{1}   &    \fbox{1} 
     \end{array}\right]$};
\end{tikzpicture}
\end{center}
\caption{\small{Analogue of Figure~\ref{Nfigure5} for matrix representations of nonobtuse $0/1$-simplices.}}
\label{Nfigure7}
\end{figure}
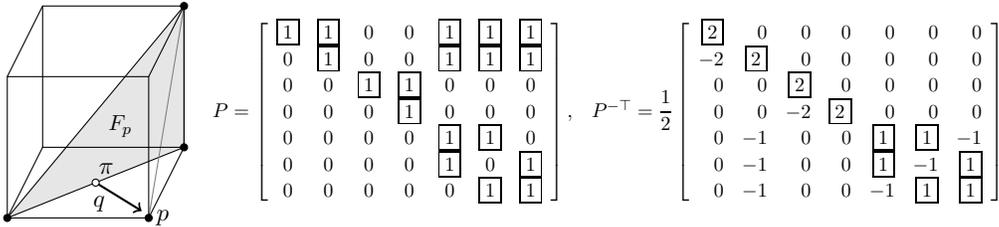   

\smallskip

Obviously, $P$ is partly decomposable and has no doubly stochastic pattern. It is only valid that the support of the doubly stochastic matrix $D$ 
is {\em contained} in the support of $P$.    
\section{Partly decomposable matrix representations}\label{Sect-4}
We will continue to study nonobtuse $0/1$-simplices having a partly decomposable matrix representation $P$. Without loss of generality, we 
may assume that $P$ is nontrivially block partitioned as
\be\label{eq-5} P = \left[\begin{array}{r|r} N & R \\ \hline 0 & A \end{array}\right], \ee
and that $A$ is fully indecomposable.

\begin{Th}\label{th-3} Let $S\in\SS^n$ be nonobtuse with a matrix representation $P$ as in (\ref{eq-5}) with $N\in\BB^{k\times k}$ with 
$k\in\{1,\dots,n-1\}$ and with $A$ fully indecomposable. Then:\\[3mm]
$\bullet$ $N$ is a matrix representation of a nonobtuse simplex in $I^k$;\\[3mm]
$\bullet$ $A$ is a matrix representation of a nonobtuse simplex in $I^{n-k}$;\\[3mm]
$\bullet$ $R=\nu(e^k)^\top$, where $\nu=0$ or $\nu$ is a column of $N$.
\end{Th}
{\bf Proof.} Lemma \ref{lem-5} proves that the first $k$ columns of $P$ together with the origin form a nonobtuse $k$-simplex, and obviously 
its vertices all lie in a $k$-facet of $I^n$. This proves the first item in the list of statements. Next, we compute 
\be\label{eq-6a} P^{-\top} = \left[\begin{array}{rc} N^{-\top} & 0 \\ -A^{-\top}R^\top N^{-\top} & A^{-\top}\end{array}\right]. \ee 
and thus,
\be\label{eq-6} (P^\top P)^{-1} = \left[\begin{array}{rr} (N^\top N)^{-1} + N^{-1}R(A^\top A)^{-1}R^\top N^{-\top} &  -N^{-1}R(A^\top A)^{-1}  \\-(A^\top A)^{-1}R^\top N^{-\top}  & (A^\top A)^{-1}\end{array}\right]. \ee 
Due to Proposition \ref{pro-1}, the matrix $(P^\top P)^{-1}$ has nonpositive off-diagonal entries (\ref{eq-1.2}) and nonnegative row sums 
(\ref{eq-1.1}). Both properties are clearly inherited by its trailing submatrix $(A^\top A)^{-1}$, possibly even with larger row sums. This 
proves the second statement of the theorem.  Next, due to (\ref{eq-1.2}), the top-right block in (\ref{eq-6}) satisfies 
\be\label{eq-7} -N^{-1}R(A^\top A)^{-1} \leq 0. \ee
Multiplication of this block from the left with $N\geq 0$ and from the right with $\ol{R}\geq 0$ gives that
\be   R(A^\top A)^{-1}\ol{R}^\top \geq 0.\ee
However, by Lemma \ref{lem-3}, the diagonal entries of $R(A^\top A)^{-1}\ol{R}^\top$ are also non{\em positive}, and thus, they all equal 
zero. We can therefore apply Lemma \ref{lem-4}. Note that because $A$ is assumed fully indecomposable, $A^\top A$ is irreducible by Lemma 
\ref{lem-1}. Thus, row-by-row application of Lemma \ref{lem-4} proves that each row of $R$ contains only zeros or only ones. This proves that 
there exists an $r\in\BB^n$ such that
\be\label{eq-8} R=r(e^{n-k})^\top.\ee
We will proceed to show that $r$ is a column of $N$, or zero. Substituting (\ref{eq-8}) back into (\ref{eq-7}) yields that  
\be\label{eq-9} wu^\top \geq 0, \hdrie \mbox{\rm where } w=N^{-1}r \und u^\top = (e^{n-k})^\top (A^\top A)^{-1}. \ee 
Due to (\ref{eq-1.1}) we have $u\geq 0$. Because $A^\top A$ is non-singular, $u$ has at least one positive entry. Thus, also $w$ is 
nonnegative. This turns $r=Nw$ into a nonnegative linear combination of columns of $N$. We continue to prove that it is a {\em convex} 
combination. For this, observe that the sums of the last $k$ rows of $(P^\top P)^{-1}$ are nonnegative due to (\ref{eq-1.1}). Thus,
\[ 0 \leq -(A^\top A)^{-1}R^\top N^{-\top}e^k + (A^\top A)^{-1}e^{n-k} = u\left(1-w^\top e^k\right) \]
with $u,w$ as in (\ref{eq-9}) and where we have used that $R^\top N^{-\top}=e^{n-k}r^\top N^{-\top} = e^{n-k}w^\top$. As we showed 
already that $u\geq 0$ has at least one positive entry, we conclude that \[  w^\top e^k \leq 1.\] 
Therefore we now have that $Nw=r\in\BB^k$ for some $w\geq 0$ with $w^\top e^k \leq 1$. Thus also
\be [0\,|\,N] \left[\begin{array}{c} 1-w^\top e^k \\ w\end{array}\right] = r. \ee
According to Lemma \ref{lem-2}, this implies that $r$ is a column of $[0\,|\,N]$. This proves the third item in the list of statements to 
prove.~\hfill $\Box$\\[3mm]
It is worthwhile to stress a number of facts concerning Theorem \ref{th-3} and its nontrivial proof.

\begin{rem}\label{rem-4}{\rm The assumption that $A$ in (\ref{eq-5}) is fully indecomposable is very natural, as each partly decomposable matrix can be put in the form (\ref{eq-5}) using operations of type (C) and (R). But in the proof of Theorem \ref{th-3} we only needed that $A^\top A$ is irreducible. This is {\em implied} by the full indecomposability of $A$ due to Lemma \ref{lem-1}, but is not {\em equivalent} to it. In fact, if $A$ is fully indecomposable, then $A^\top A\geq e^n(e^n)^\top + I$. See Corollary \ref{co-5}.}
\end{rem} 

\begin{rem}\label{rem-5}{\rm The result proved in the third bullet of Theorem \ref{th-3} that $R$ consists of $n-k$ copies of {\em the same} column of $N$ is stronger than the geometrical observation that each of the last $n-k$ columns of $P$ should project on {\em any} vertex of the $k$-simplex represented by $N$. It is the {\em irreducibility} of $A^\top A$ that forces the equality of all columns of $R$.}
\end{rem}

\begin{rem}\label{rem-6}{\rm Permuting rows and columns of the block upper triangular matrix in (\ref{eq-5}) show that also for
\be  \left[\begin{array}{r|r} A & 0 \\ \hline R & N \end{array}\right]  \ee
with $A$ and $N$ as in Theorem \ref{th-3}, similar conclusions can be drawn for $R$}.
\end{rem}
Some details in the above remarks will turn out to be of central importance in Section \ref{Sect-6}.

\begin{Co}\label{co-3} Let $S\in\SS^n$ be a nonobtuse $0/1$-simplex with matrix representation $P$. Then the following statements are 
equivalent:\\[3mm]
$\bullet$ $P$ is partly decomposable;\\[3mm]
$\bullet$ $S$ has a block diagonal matrix representation with at least one fully indecomposable block;\\[3mm]
$\bullet$ each matrix representation of $S$ is partly decomposable.
\end{Co} 
{\bf Proof}. Suppose that $P$ is partly decomposable. Then Theorem \ref{th-1} shows that $P$ is of the form
\be\label{eq-10} P = \left[\begin{array}{r|c} N & \nu \left(e^{n-k}\right)^\top \\[1mm] \hline 0 & A \end{array}\right],\ee
and $\nu=0$ or $\nu$ is a column of $N$. If $\nu=0$ then $P$ itself is block diagonal. If $\nu\not=0$, apply to $P$ the operation of type (X) 
as described below Figure~\ref{Nfigure3} with column $c$ equal to the column $(\nu,0)^\top$ of $P$. The simple observation that $\nu+\nu$ equal zero 
modulo $2$ proves that the resulting matrix $\tilde{P}$ is block diagonal. As the bottom right block of $\tilde{P}$ equals $A$, this shows that 
at least one block is fully indecomposable. To show that each matrix representation of $S$ is partly decomposable, simply note that each 
operation of type (X) applied to the block-diagonal matrix representation will leave one of the two off-diagonal zero blocks 
invariant.~\hfill $\Box$

\begin{rem}{\rm The converse of Theorem \ref{th-3} is also valid. Indeed, suppose that $N$ and $A$ are matrix representations of nonobtuse 
simplices. Then it is trivially true that the block diagonal matrix $P$ having $N$ and $A$ as diagonal blocks represents a nonobtuse simplex. 
Applying operations of type (X) to $P$ proves that all matrices of the form (\ref{eq-10}) then represent nonobtuse simplices. Note that this also 
holds without the assumption that $A$ is fully indecomposable.}
\end{rem}
Another corollary of Theorem \ref{th-3} concerns its implications for the structure of the transposed inverse $P^{-\top}$ of a partly 
decomposable matrix representation of a nonobtuse $0/1$-simplex.

\begin{figure}[h]
\begin{center}
\begin{tikzpicture}[scale=0.94, every node/.style={scale=0.94}]

\node[xscale=0.6,yscale=0.73] at (5,1)  {$\left[\begin{array}{ccccccc}
     \blacksquare  &   \blacksquare   &   \square &  \square &    \blacksquare   &    \blacksquare   &    \blacksquare \\
    \square &    \blacksquare   &  \square &  \square  &    \blacksquare   &    \blacksquare   &    \blacksquare\\
    \square &  \square &    \blacksquare   &   \blacksquare   &   \square  &  \square  &   \square\\
    \square &  \square &  \square &    \blacksquare   &  \square &  \square &  \square\\
    \square &  \square &  \square &  \square &    \blacksquare   &    \blacksquare   &   \square\\ 
    \square &  \square &  \square &  \square &    \blacksquare   &  \square &    \blacksquare \\
    \square &  \square &  \square &  \square &  \square &    \blacksquare   &    \blacksquare 
     \end{array}\right]$};
     
\draw[<->] (6.8,1)--(7.8,1);
\draw[<->] (3.3,1)--(2.3,1);
\node[scale=0.9] at (7.3,1.3) {(X)};
\node[scale=0.9] at (2.8,1.3) {(X)};
\node[scale=0.9] at (7.3,0.7) {$c=2$};
\node[scale=0.9] at (2.8,0.7) {$c=6$};

\draw[<->] (10,-1.3)--(9,-2.3);
\draw[<->] (0,-1.3)--(1,-2.3);
\node[scale=0.9] at (10.5,-2) {(R), $1\leftrightarrow 2$};
\node[scale=0.9] at (-0.5,-2) {(C), $2\leftrightarrow 6$};
\node[scale=0.9] at (-0.5,-2.5) {(R), $1\leftrightarrow 2$};

\draw[<->] (4.5,-3)--(5.5,-3);
\node[scale=0.9] at (5,-2.7) {(X)};
\node[scale=0.9] at (5,-3.3) {$c=6$};
     
\node[xscale=0.6,yscale=0.73] at (9.5,0)  {$\left[\begin{array}{ccccccc}
    \square&    \blacksquare   &   \blacksquare  &   \blacksquare   &  \square  &   \square  &  \square\\
       \blacksquare &   \blacksquare  &    \blacksquare  &   \blacksquare  &  \square  &  \square  &  \square\\
    \square &  \square &    \blacksquare   &   \blacksquare   &   \square  &  \square  &   \square\\
    \square &  \square &  \square &    \blacksquare   &  \square &  \square &   \square\\
    \square &  \square &  \square &  \square &    \blacksquare   &    \blacksquare   &   \square\\ 
    \square &  \square &  \square &  \square &    \blacksquare   &  \square &    \blacksquare \\
    \square &  \square &  \square &  \square &  \square &    \blacksquare   &    \blacksquare 
     \end{array}\right]$};
     
 \node[xscale=0.6,yscale=0.73] at (7.2,-3)  {$\left[\begin{array}{ccccccc}
     \blacksquare &   \blacksquare  &    \blacksquare  &   \blacksquare  &  \square  &  \square  &  \square\\
    \square&    \blacksquare   &   \blacksquare  &   \blacksquare   &  \square  &   \square  &  \square\\
    \square &  \square &    \blacksquare   &   \blacksquare   &   \square  &  \square  &   \square\\
    \square &  \square &  \square &    \blacksquare   &  \square &  \square & \square\\
    \square &  \square &  \square &  \square &    \blacksquare   &    \blacksquare   &  \square\\ 
    \square &  \square &  \square &  \square &    \blacksquare   &  \square &    \blacksquare \\
    \square &  \square &  \square &  \square &  \square &    \blacksquare   &    \blacksquare 
     \end{array}\right]$};

\node[xscale=0.6,yscale=0.73] at (0.5,0)  {$\left[\begin{array}{ccccccc}
     \square  &   \square   &   \blacksquare &  \blacksquare &   \square  &   \blacksquare &  \square\\
    \blacksquare &    \square   &  \blacksquare &  \blacksquare  &   \square &  \blacksquare  &  \square\\
    \square &  \square &    \blacksquare   &   \blacksquare   &   \square  &  \square  &   \square\\
    \square &  \square &  \square &    \blacksquare   &  \square &  \square & \square\\
     \blacksquare  &   \blacksquare  &   \blacksquare  &   \blacksquare  &   \square  &    \blacksquare   &   \blacksquare\\ 
    \square &  \square &  \square &  \square &    \blacksquare   &  \square &    \blacksquare \\
     \blacksquare &   \blacksquare  &   \blacksquare  &   \blacksquare  &   \blacksquare  &    \blacksquare  &   \square
     \end{array}\right]$};
     
\node[xscale=0.6,yscale=0.73] at (2.8,-3)  {$\left[\begin{array}{ccccccc}
    \blacksquare &    \blacksquare   &  \blacksquare &  \blacksquare  &   \square &  \square  &  \square\\
       \square  &   \blacksquare   &   \blacksquare &  \blacksquare &   \square  &   \square &  \square\\
    \square &  \square &    \blacksquare   &   \blacksquare   &   \square  &  \square  &   \square\\
    \square &  \square &  \square &    \blacksquare   &  \square &  \square & \square\\
     \blacksquare  &   \blacksquare  &   \blacksquare  &   \blacksquare  &   \square  &    \blacksquare   &   \blacksquare\\ 
    \square &  \square &  \square &  \square &    \blacksquare   &  \square &    \blacksquare \\
     \blacksquare &   \blacksquare  &   \blacksquare  &   \blacksquare  &   \blacksquare  &    \blacksquare  &   \square
     \end{array}\right]$};

\end{tikzpicture}
\end{center}
\caption{\small{Illustration of Corollary \ref{co-3} using the matrix $P\in\BB^{7\times7}$ from Figure~\ref{Nfigure7}.}}
\label{Nfigure8}
\end{figure}
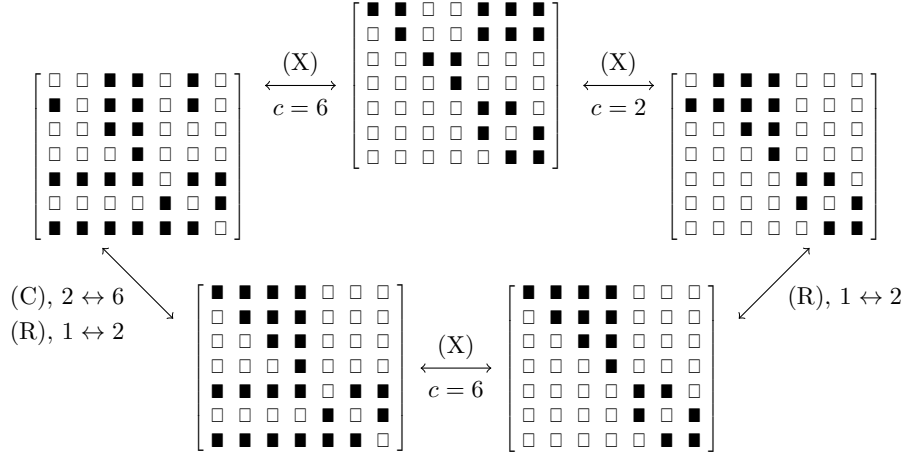

\begin{Co}\label{co-4} If $\nu=Ne_j^k$ in (\ref{eq-10}) for some $j\in\{1,\dots,k\}$, then
\be P^{-\top} = \left[\begin{array}{c|c} N^{-\top} & 0 \\\hline\\[-4mm] ae_j^\top & A^{-\top}\end{array}\right],  \ee
where $a$ is the inward normal to the facet opposite the origin of the simplex represented by $A$. As a consequence, the sums of the 
last $n-k$ rows of $P^{-\top}$ all add to zero.
\end{Co}
{\bf Proof. } Substitute $R=\nu e^{n-k}$ with $\nu=Ne_j^k$ into the expression for $P^{-\top}$ in (\ref{eq-6a}). \hfill $\Box$\\[3mm]
To illustrate Corollary \ref{co-3}, consider again the matrix $P\in\BB^{7\times7}$ from Figure~\ref{Nfigure7}, now displayed not as $0/1$-matrix but as 
checkerboard black-white pattern, at the top in Figure~\ref{Nfigure8}. Applying an operation of type (X) with the second column results in the matrix to its 
right, which is block diagonal. Swapping the first two rows of that matrix yields one in which the top left block is now in its most reduced form. 
Applying an operation of type (X) with the sixth column results in a matrix in which the bottom left $3\times 4$ zero block has been destroyed. 
Finally, swapping columns $2$ and $6$, and swapping rows $1$ and $2$, results in the matrix that could also have been obtained by applying 
operation (X) with column $6$ directly to $P$.

\smallskip

Corollary \ref{co-3} shows that the matrix representations of a nonobtuse $0/1$-simplex are either all partly decomposable, or all fully 
indecomposable. This motivates to the following definition.
 
\begin{Def} {\rm A nonobtuse simplex is called partly decomposable if it has a partly decomposable matrix representation, and fully 
indecomposable if it has not}.
\end{Def}
We will now investigate to what structure the recursive application of Theorem \ref{th-3} leads. For this, assume again that $P$ is a partly 
decomposable matrix representation of a nonobtuse $0/1$-simplex $S\in\SS^n$. Then by Corollary \ref{co-3}, $S$ has a matrix representation 
of the form
\be\label{eq-11} P = \left[\begin{array}{r|c} N_1 & R_1 \\ \hline 0 & A_1 \end{array}\right],\ee
in which $A_1$ is fully indecomposable. According to Theorem \ref{th-3}, the $k\times k$ matrix $N_1$ represents a nonobtuse $k$-simplex 
$K$ in $I^k$. If also $K$ is partly decomposable, we can block-partition $N_1$ using row (R) and column (C) permutations $\Pi_1$ and 
$\Pi_2$, such that 
\be\label{eq-12} \tilde{P}=\Pi_1 P\Pi_2 = \left[\begin{array}{c|c|c} N_2 & R_{12} & R_{13} \\\hline 0 & A_2 & R_{23}\\ \hline 0 & 0 & A_1 \end{array}\right],\ee
with $A_2$ fully indecomposable and $N_2$ possibly partly decomposable. Theorem \ref{th-3} shows that
\[ R_{12} \hdrie \mbox{\rm consists of copies of a column $\nu$ of} \hdrie [0|N_2], \]
and
\[ \left[\begin{array}{c}R_{13}\\\hline R_{23}\end{array}\right]\hdrie\mbox{\rm consists of copies of a column of} \hdrie \left[\begin{array}{r|c} N_2 & R_{12} \\ \hline 0 & A_2 \end{array}\right]. \]
This means that if $R_{23}$ is nonzero, then $R_{13}$ consists of copies of the same column $\nu$ of $N$ as does $R_{12}$. Thus the whole 
block $[R_{12}\,|\,R_{23}]$ consists of copies of a column $\nu$ of $[0|N_2]$. On the other hand, if $R_{23}$ is zero, then $R_{13}$ can 
either be zero, or consist of copies of {\em any} column of $N_2$, including $\nu$.\\[3mm]
By including operations of type (X) it is possible to map the entire strip above one of the fully indecomposable diagonal blocks to zero. Although 
this will in general destroy the block upper triangular form of the square submatrix to the left of that strip, Corollary \ref{co-4} shows that this 
submatrix remains partly decomposable. Therefore, its block upper triangular structure can be restored using only operations of type (R) and 
(C), which leave the zero strip intact.

\begin{rem}{\rm It is generally not possible to transform $\tilde{P}$ to block diagonal form with more than two diagonal blocks. See the 
tetrahedron $S$ in $I^3$ in Figure~\ref{Nfigure9}. It is possible to put a vertex of $S$ at the origin such that facets $A$ and $N$ are orthogonal, and 
hence the corresponding matrix representation is block diagonal with two blocks. The triangular facet however requires a {\em different} vertex 
at the origin for its $2\times 2$ matrix representation to be diagonal.}
\end{rem}
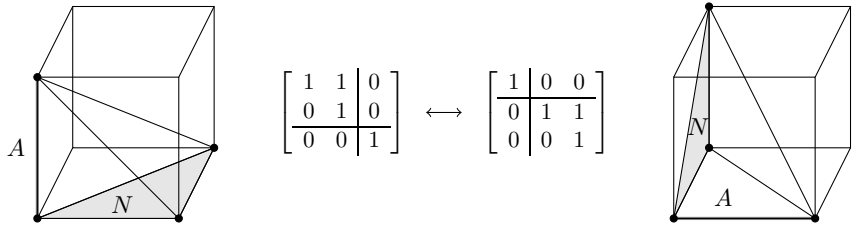
\begin{figure}[h]
\begin{center}
\begin{tikzpicture}[scale=0.94, every node/.style={scale=0.94}] 
\begin{scope}[shift={(0,0)}]
\draw[fill=gray!20!white] (0,0)--(2,0)--(2.5,1)--(0,0)--cycle;
\draw (2,0)--(0,2);
\draw (0,0)--(2.5,1);
\draw (0,2)--(2.5,1);
\draw[gray,very thick] (0,0)--(0,2);
\draw (0,0)--(2,0)--(2,2)--(0,2)--cycle;
\draw (0.5,1)--(2.5,1)--(2.5,3)--(0.5,3)--cycle;
\draw (0,0)--(0.5,1);
\draw (2,0)--(2.5,1);
\draw (0,2)--(0.5,3);
\draw (2,2)--(2.5,3);
\draw[fill=black] (0,0) circle [radius=0.05];
\draw[fill=black] (2,0) circle [radius=0.05]; 
\draw[fill=black] (2.5,1) circle [radius=0.05];
\draw[fill=black] (0,2) circle [radius=0.05];
\node[scale=0.9] at (-0.3,1) {$A$};
\node[scale=0.9] at (1.2,0.2) {$N$};
\end{scope}
\node[scale=0.85] at (4.3,1.5) {$\left[\begin{array}{rr|r} 1 & 1 & 0 \\ 0 & 1 & 0 \\\hline 0 & 0 & 1\end{array}\right]$};
\draw[<->] (5.5,1.5)--(6,1.5);
\node[scale=0.85] at (7.2,1.5) {$\left[\begin{array}{r|rr} 1 & 0 & 0 \\\hline 0 & 1 & 1 \\ 0 & 0 & 1\end{array}\right]$};
\begin{scope}[shift={(9,0)}]
\draw[fill=gray!20!white] (0,0)--(0.5,1)--(0.5,3)--cycle;
\draw (2,0)--(0.5,1);
\draw (2,0)--(0.5,3);
\draw[gray,very thick] (0,0)--(2,0);
\draw (0,0)--(2,0)--(2,2)--(0,2)--cycle;
\draw (0.5,1)--(2.5,1)--(2.5,3)--(0.5,3)--cycle;
\draw (0,0)--(0.5,1);
\draw (2,0)--(2.5,1);
\draw (0,2)--(0.5,3);
\draw (2,2)--(2.5,3);
\draw[fill=black] (0,0) circle [radius=0.05];
\draw[fill=black] (2,0) circle [radius=0.05]; 
\draw[fill=black] (0.5,3) circle [radius=0.05];
\draw[fill=black] (0.5,1) circle [radius=0.05];
\node[scale=0.9] at (0.7,0.3) {$A$};
\node[scale=0.9] at (0.35,1.3) {$N$};
\end{scope}
\end{tikzpicture}
\end{center}
\caption{\small{Any simplex $S$ with a partly decomposable matrix representation has a pair of facets $A$ and $N$ of dimensions adding to 
$n$ are orthogonal to one another.}}
\label{Nfigure9}
\end{figure}      

\smallskip

Summarizing, the above discussion shows that each nonobtuse $0/1$-simplex $S$ has a matrix representation that is block upper triangular, 
with fully indecomposable diagonal blocks (possibly only one). The strip above each diagonal block is of rank one and consists only of copies of 
a column to the left of the strip. Any matrix representation $P$ of $S$ can be brought into this form using only operations of type (C) and (R). 
Using an additional reflection of type (X), it is possible to transform an entire strip above one of the diagonal blocks to zero using a column to 
the left of the strip. Although this may destroy the block upper triangular form to the left of the strip, this form can be restored using operations 
of type (R) and (C) only.\\[3mm]
In view of Corollary \ref{co-4}, the transposed inverse $P^{-\top}$ of a partly decomposable matrix representation $P$ of a nonobtuse $0/1$-
simplex $S$ is perhaps even simpler in structure than $P$ itself. On the diagonal it has the transposed inverses of the fully indecomposable 
diagonal blocks $A_1,\dots,A_p$ of $P$, and each {\em horizontal} strip to the left of such a diagonal block $(A_j)^{-\top}$ has 
{\em at most} one nonpositive column that is not identically zero. This columns has two interesting features. The first is that it nullifies the 
sums of the 
rows in its strip. The second is that its position clearly indicates to which vertex of which other block $A_i$ the block $A_j$ is related. To 
illustrate what we mean by this, see Figure~\ref{Nfigure7} for an example. Above the bottom right $3\times 3$ block $A$ of $P$ we see copies of the second 
column of $P$. This fact can also be read from the position of the nonzero entries to the left of the corresponding block $A^{-\top}$ in 
$P^{-\top}$.
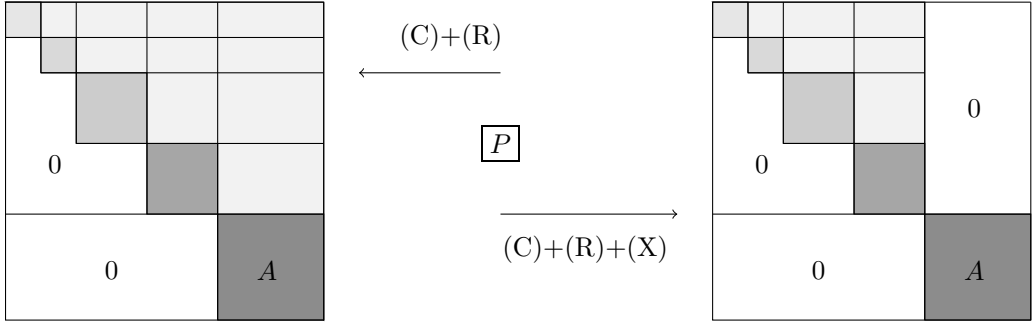
\begin{figure}[h]
\begin{center}
\begin{tikzpicture}[scale=0.94, every node/.style={scale=0.94}]
\draw[fill=gray!10!white] (0,4.5)--(4.5,0)--(4.5,4.5)--cycle;
\draw[fill=gray!90!white] (3,0)--(4.5,0)--(4.5,1.5)--(3,1.5)--cycle;
\draw[fill=gray!70!white] (2,1.5)--(3,1.5)--(3,2.5)--(2,2.5)--cycle;
\draw[fill=gray!40!white] (1,2.5)--(2,2.5)--(2,3.5)--(1,3.5)--cycle;
\draw[fill=gray!30!white] (0.5,3.5)--(1,3.5)--(1,4)--(0.5,4)--cycle;
\draw[fill=gray!20!white] (0,4)--(0.5,4)--(0.5,4.5)--(0,4.5)--cycle;
\node at (0.7,2.2) {$0$};
\node at (1.5,0.7) {$0$};
\node at (3.7,0.7) {$A$}; 
\draw (0,4)--(0,0)--(3,0);
\draw (3,0)--(4.5,0);
\draw (0,1.5)--(4.5,1.5);
\draw (1,2.5)--(4.5,2.5);
\draw (0.5,3.5)--(4.5,3.5);
\draw (0,4)--(4.5,4);
\draw (0,4.5)--(4.5,4.5);
\draw (0,4)--(0,4.5);
\draw (0.5,3.5)--(0.5,4.5);
\draw (1,2.5)--(1,4.5);
\draw (2,1.5)--(2,4.5);
\draw (3,0)--(3,4.5);
\draw (4.5,0)--(4.5,4.5);
\node at (7,2.5) {\fbox{$P$}};
\draw[->] (7,3.5)--(5,3.5);
\draw[->] (7,1.5)--(9.5,1.5);
\node at (6.3,4) {(C)+(R)};
\node at (8.2,1) {(C)+(R)+(X)};
\begin{scope}[shift={(10,0)}]
\draw[fill=gray!10!white] (0,4.5)--(3,1.5)--(3,4.5)--cycle;
\draw[fill=gray!90!white] (3,0)--(4.5,0)--(4.5,1.5)--(3,1.5)--cycle;
\draw[fill=gray!70!white] (2,1.5)--(3,1.5)--(3,2.5)--(2,2.5)--cycle;
\draw[fill=gray!40!white] (1,2.5)--(2,2.5)--(2,3.5)--(1,3.5)--cycle;
\draw[fill=gray!30!white] (0.5,3.5)--(1,3.5)--(1,4)--(0.5,4)--cycle;
\draw[fill=gray!20!white] (0,4)--(0.5,4)--(0.5,4.5)--(0,4.5)--cycle;
\node at (1.5,0.7) {$0$};
\node at (3.7,3) {$0$};
\node at (3.7,0.7) {$A$};
\node at (0.7,2.2) {$0$};
\draw (0,4)--(0,0)--(3,0);
\draw (3,0)--(4.5,0);
\draw (0,1.5)--(4.5,1.5);
\draw (1,2.5)--(3,2.5);
\draw (0.5,3.5)--(3,3.5);
\draw (0,4)--(3,4);
\draw (0,4.5)--(4.5,4.5);
\draw (0,4)--(0,4.5);
\draw (0.5,3.5)--(0.5,4.5);
\draw (1,2.5)--(1,4.5);
\draw (2,1.5)--(2,4.5);
\draw (3,0)--(3,4.5);
\draw (4.5,0)--(4.5,4.5);
\end{scope}

\end{tikzpicture}
\end{center}
\caption{\small{A partly decomposable matrix representation $P$ of a nonobtuse $0/1$-simplex $S$ can be brought in the left form using 
operations of type (C) and (R) only, and in the right form if using additional operations of type (X). The diagonal blocks are fully 
indecomposable.}}
\label{Nfigure10}
\end{figure}  

\begin{rem}\label{rem-10}{\rm The above shows that to each nonobtuse $0/1$-simplex $S$ of dimension $n$ we can associate a special type 
of simplicial complex $C_p$ consisting of $p$ mutually orthogonal fully indecomposable simplicial facets $S_1,\dots, S_p$ with respective 
dimensions $k_1,\dots,k_p$ adding to $n$, where each facet $S_j$ lies in in its own $k_j$-facet of $I^n$. Explicitly, let $C_1=S_1$. The 
complex $C_{j+1}$ is obtained by attaching a vertex of $S_{j+1}$ to a vertex $v$ of $C_j$, such that the orthogonal projection of $S_{j+1}$ 
onto the $(k_1+\dots+k_j)$-dimensional ambient space of $C_j$ equals $v$.}
\end{rem}
Remark \ref{rem-10} is illustrated by the tetrahedron in Figure~\ref{Nfigure9}. It can be built from three $1$-simplices $S_1,S_2,S_3$ simply by first 
attaching $S_2$ with a vertex to a vertex of $S_1$ orthogonally to $S_1$, giving a right triangle $C_2$. Then attaching $S_3$ to the correct 
vertex $v$ of $C_2$ such that the projection of $S_3$ onto $C_2$ equals $v$ gives the tetrahedron. In Section~\ref{Sect-5} we will pay special 
attention to the nonobtuse simplices whose fully indecomposable components are $n$-cube edges.

\begin{figure}[h]
\begin{center}
\begin{tikzpicture}[scale=0.94, every node/.style={scale=0.94}]
\begin{scope}[scale=1.1, every node/.style={scale=1.1}]
\draw[fill=gray!30!white] (2,2)--(3.8,0.5)--(5,1.3)--(5,2.7)--(3.8,3.5)--cycle;
\draw[fill=gray!30!white] (1,2)--(-1,1)--(-1.5,2.5)--(-0.2,3.4)--cycle;
\draw (2,2)--(1,2);
\node[scale=0.9] at (3.2,2) {$S_1$};
\node[scale=0.9] at (-0.3,1.9) {$S_3$};
\node[scale=0.9] at (1.5,1.7) {$S_2$};
\draw (3.8,0.5)--(5,2.7)--(2,2); 
\draw (2,2)--(5,1.3)--(3.8,3.5); 
\draw (3.8,0.5)--(3.8,3.5);
\draw (-1,1)--(-0.2,3.4);
\draw[gray] (1,2)--(-1.5,2.5);
\draw[fill=black] (-1,1) circle [radius=0.05];
\draw[fill=white] (1,2) circle [radius=0.05];
\draw[fill=black] (2,2) circle [radius=0.05];
\draw[fill=black] (-1.5,2.5) circle [radius=0.05]; 
\draw[fill=black] (-0.2,3.4) circle [radius=0.05];
\draw[fill=black] (3.8,0.5) circle [radius=0.05];
\draw[fill=black] (5,1.3) circle [radius=0.05];
\draw[fill=black] (5,2.7) circle [radius=0.05];
\draw[fill=black] (3.8,3.5) circle [radius=0.05];
\node[scale=0.8] at (0,3.6) {$1$};
\node[scale=0.8] at (1.2,2.2) {$0$};
\node[scale=0.8] at (1.8,2.2) {$4$};
\node[scale=0.8] at (-1,0.7) {$3$};
\node[scale=0.8] at (-1.7,2.5) {$2$};
\node[scale=0.8] at (4,0.4) {$5$};
\node[scale=0.8] at (5.2,1.4) {$6$};
\node[scale=0.8] at (5.2,2.6) {$7$};
\node[scale=0.8] at (4,3.6) {$8$};
\end{scope}
\begin{scope}[shift={(8.5,2.3)}]
\node[scale=0.9] at (0,-1.5) {$1\,\,\,\,2\,\,\,\,3\,\,\,\,4\,\,\,\,5\,\,\,\,6\,\,\,\,7\,\,\,8$}; 
\node[xscale=0.6,yscale=0.73] at (0,0.1)  {$\left[\begin{array}{ccc|c|cccc}
     \blacksquare  &   \blacksquare   &   \square &  \square &    \square   &    \square   &    \square &    \square\\
    \blacksquare &    \square   &  \blacksquare &  \square  &    \square   &    \square   &    \square&    \square\\
    \square &  \blacksquare &    \blacksquare   &   \square   &   \square  &  \square  &   \square&    \square\\\hline
    \square &  \square &  \square &    \blacksquare   &  \blacksquare &  \blacksquare &  \blacksquare&    \blacksquare\\\hline
    \square &  \square &  \square &  \square &    \blacksquare    &   \blacksquare   &    \blacksquare &   \square\\ 
    \square &  \square &  \square &  \square  &   \blacksquare     &  \blacksquare         &     \square     &  \blacksquare \\
    \square &  \square &  \square &  \square  &   \blacksquare          &  \square    &    \blacksquare      & \blacksquare \\
     \square &  \square &  \square &  \square &   \square          &  \blacksquare         & \blacksquare    & \blacksquare \\
     \end{array}\right]$};      
     \node[scale=0.9] at (3.3,0.1) {$=\left[\begin{array}{c|c|c} A_3 & 0 & 0 \\\hline 0 & A_2 & R \\\hline 0 & 0 & A_1\end{array}\right]$};
\end{scope}
\end{tikzpicture}
\end{center}
\caption{\small{A simplicial complex of three mutually orthogonal simplices $S_1,S_2,S_3$. The given matrix representation corresponds to 
choosing the vertex $0$ as the origin. Reflecting vertex $4$ to the origin decouples, alternatively, the bottom right $4\times 4$ block. It is not 
possible to transform the matrix to block diagonal form with $A_2=[\,1\,]$ as one of the diagonal blocks. Reflecting any other vertex to the 
origin does not even lead to a block diagonal matrix representation.}}
\label{Nfigure11}
\end{figure}
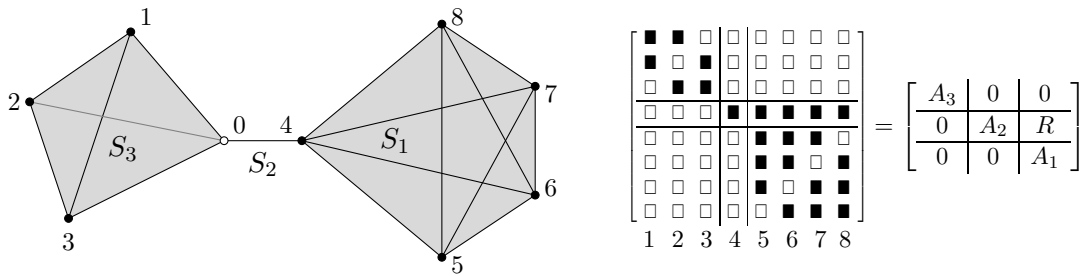

\begin{rem}{\rm The $1$-simplex in $I^n$ has a fully indecomposable matrix representation with doubly stochastic pattern. It is formally an 
acute simplex. Indeed, the normals to its $0$-dimensional facets $0$ and $1$ point in opposite directions. Hence, its only dihedral angle equals 
zero. Since there does not exist a fully indecomposable triangle in $I^2$, matrix representations of a partly decomposable simplex do not have 
$2\times 2$ fully indecomposable diagonal blocks.}
\end{rem}  
We would like to stress that although each partly decomposable matrix representation of a nonobtuse $0/1$-simplex can be transformed into 
block diagonal form by operations of types (C),(R) and (X) as depicted in the right in Figure~\ref{Nfigure10}, the bottom right block cannot be {\em any} of 
the fully indecomposable diagonal blocks $A_j$. This is only possible if the corresponding simplex $S_j$ is attached to the remainder of the 
complex at exactly one vertex. For example, in Figure~\ref{Nfigure11}, with $0$ as the origin, the matrix $A_3$ representing a regular tetrahedron in $I^3$ 
is a block of a block diagonal matrix. After mapping vertex $4$ to the origin with a reflection of type (X), the block $A_4$ representing a so-
called antipodal $4$-simplex in $I^4$ is. However, the block $A_2$ representing the $1$-simplex in $I^1$ is never a block of a block diagonal 
matrix representation of $S$. The only configurations of the three building blocks $S_1,S_2,S_3$ having a matrix representation that can be 
transformed by operations of type (C),(R) and (X) onto block diagonal form with three diagonal blocks, are those in which $S_1,S_2$ and 
$S_3$ have a common vertex.

\smallskip

\begin{figure}[h]
\begin{center}
\begin{tikzpicture}[scale=0.94, every node/.style={scale=0.94}]
\draw[fill=gray!10!white] (0,0)--(1,-1)--(2,0)--(1,1)--cycle;
\draw (2,0)--(3,0);
\draw[fill=gray!30!white] (3,0)--(2,-0.7)--(2.4,-1.7)--(3.6,-1.7)--(4,-0.7)--cycle;
\draw (3,0)--(4,0);
\draw (4,0)--(5,1);
\draw (4,0)--(3,1);
\draw[fill=gray!50!white] (5,1)--(5.5,1.7)--(6.2,1.7)--(6.7,1)--(6.2,0.3)--(5.5,0.3)--cycle;
\draw[fill=gray!20!white] (5,0.3)--(6,-0.7)--(5,-1.7)--(4,-0.7)--cycle;
\draw[fill=black] (0,0) circle [radius=0.05];
\draw[fill=white] (2,0) circle [radius=0.05];
\draw[fill=white] (3,0) circle [radius=0.05];
\draw[fill=white] (4,0) circle [radius=0.05];
\draw[fill=white] (5,1) circle [radius=0.05];
\draw[fill=black] (3,1) circle [radius=0.05];
\draw[fill=black] (1,1) circle [radius=0.05];
\draw[fill=black] (1,0-1) circle [radius=0.05];
\draw[fill=black] (2,-0.7) circle [radius=0.05];
\draw[fill=black] (2.4,-1.7) circle [radius=0.05];
\draw[fill=black] (3.6,-1.7) circle [radius=0.05];
\draw[fill=black] (4,-0.7) circle [radius=0.05];
\draw[fill=black] (5.5,1.7) circle [radius=0.05];
\draw[fill=black] (6.2,1.7) circle [radius=0.05];
\draw[fill=black] (6.2,0.3) circle [radius=0.05];
\draw[fill=black] (5.5,0.3) circle [radius=0.05];
\draw[fill=black] (6.7,1) circle [radius=0.05];
\draw[fill=black] (6,-0.7) circle [radius=0.05];
\draw[fill=black] (5,-1.7) circle [radius=0.05];
\draw[fill=black] (5,0.3) circle [radius=0.05];
\draw[fill=white] (4,-0.7) circle [radius=0.05];
\node[scale=0.8] at (2.1,0.2) {$2$};
\node[scale=0.8] at (3,0.2) {$3$};
\node[scale=0.8] at (4,0.3) {$3$};
\node[scale=0.8] at (4.7,1) {$2$};
\node[scale=0.8] at (4,-0.4) {$2$};
\end{tikzpicture}
\end{center}
\caption{\small{Schematic representation of a simplicial complex, built from a $5$-simplex, a $4$-simplex, two tetrahedra, and four edges, of 
total dimension $19$. With the origin at a white vertex, the matrix representation decouples into the indicated number of diagonal blocks. If the 
vertex is located at another vertex, the matrix representation does not decouple.}}
\label{Nfigure12}
\end{figure}
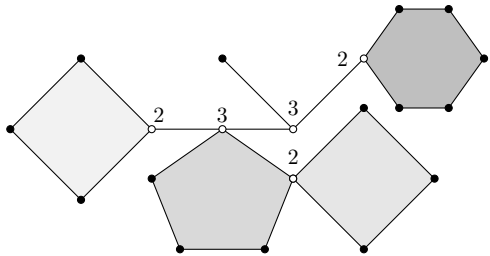

In general, the decomposability structure of a nonobtuse $0/1$-simplex can be well visualized as a special type of planar graph, at the cost of 
the geometrical structure. For this, assign to each $p\times p$ fully indecomposable diagonal block a regular $p$-gon, and to attach these to 
one another at the common vertex of the simplices they represent. See Figure~\ref{Nfigure12} for an example. At the white vertices it is indicated how many 
$p$-gons meet. This number equals the number of diagonal blocks in the matrix representation when this vertex is reflected onto the origin. 

\smallskip

Before studying further properties of partly decomposable nonobtuse $0/1$-simplices in terms of their fully indecomposable components, we 
will pay special attention to {\em orthogonal} simplices.
\section{Orthogonal simplices and their matrix representations}\label{Sect-5}
The simplest class of nonobtuse simplices is formed by the {\em orthogonal simplices}. These are nonobtuse simplices with $\binom{n}{2}-n$ 
right dihedral angles. Note that this is the {\em maximum} number of right dihedral angles a simplex can have, as Fiedler proved in \cite{Fie} 
that any simplex has at least $n$ acute dihedral angles. Orthogonal simplices are useful in many applications, see \cite{BrDiHaKr,BrKoKrSo} 
and the references therein. We will restrict our attention to orthogonal $0/1$-simplices.\\[3mm]
The orthogonal $0/1$-simplices can be defined recursively as follows \cite{BrDiHaKr}. The cube edge $I^1$ is an orthogonal simplex. Now, a 
nonobtuse $0/1$-simplex $S$ in $I^n$ is orthogonal if it has an $(\nmo)$-facet $F$ with the properties that:\\[2mm]
$\bullet$ $F$ is contained in an $(\nmo)$-facet of $I^n$;\\[2mm]
$\bullet$ $F$ is an orthogonal $(\nmo)$-simplex.\\[2mm]
Clearly, this way to construct orthogonal $0/1$-simplices is a special case of how nonobtuse $0/1$-simplices were constructed from their fully 
indecomposable parts in Section~\ref{Sect-4}. This is because a vertex $v$ forms an orthogonal simplex $S$ together with an $(\nmo)$-facet 
$F$ that is contained in an $(\nmo)$-facet of $I^n$ if and only if $v$ projects orthogonally on a vertex of $F$. This shows in particular that all 
fully indecomposable components of any matrix representation of $S$ equal $[\,1\,]$, which limits the number of their upper triangular matrix 
representations.

\begin{Pro} There exist $n!$ distinct upper triangular $0/1$ matrices that represent orthogonal $0/1$-simplices in $I^n$. 
\end{Pro}
{\bf Proof. } Let $P\in\BB^{n\times n}$ be an upper triangular matrix representing an orthogonal $n$-simplex. Then the matrix
\[ \tilde{P} = \left[\begin{array}{r|r} P & r \\\hline 0 & 1\end{array}\right] \]
is upper triangular, and according to Theorem \ref{th-3} it represents a nonobtuse simplex if and only if $r$ equals one of the $n+1$ distinct 
columns of $[0\,|\,P]$. Thus there are $\npo$ times as many matrix representations of orthogonal $(\npo)$-simplices in $I^{\npo}$ as of 
orthogonal $n$-simplices in $I^n$, whereas $[\,1\,]$ is the only one in $I^1$.\hfill $\Box$
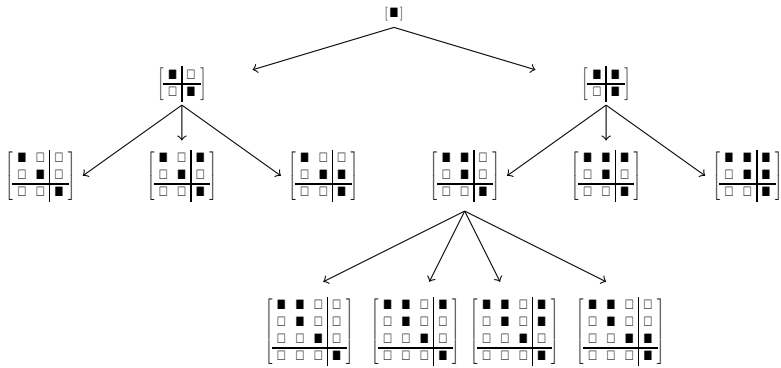
\begin{figure}[h]
\begin{center}
\begin{tikzpicture}[scale=0.94, every node/.style={scale=0.94}]
\draw[->] (6,-0.2)--(4,-0.8);
\draw[->] (6,-0.2)--(8,-0.8);
\draw[->] (3,-1.3)--(1.6,-2.3);
\draw[->] (3,-1.3)--(4.4,-2.3);
\draw[->] (3,-1.3)--(3,-1.8);
\draw[->] (9,-1.3)--(7.6,-2.3);
\draw[->] (9,-1.3)--(10.4,-2.3);
\draw[->] (9,-1.3)--(9,-1.8);
\draw[->] (7,-2.8)--(5,-3.8);
\draw[->] (7,-2.8)--(6.5,-3.8);
\draw[->] (7,-2.8)--(7.5,-3.8);
\draw[->] (7,-2.8)--(9,-3.8);
\node[xscale=0.4, yscale=0.5]  at (6,0) {$[\,\blacksquare\,]$};
\node[xscale=0.4, yscale=0.5]  at (3,-1) {$\left[\begin{array}{c|c} \blacksquare & \square \\ \hline \square & \blacksquare\end{array}\right]$};
\node[xscale=0.4, yscale=0.5]  at (9,-1) {$\left[\begin{array}{c|c} \blacksquare & \blacksquare \\ \hline \square & \blacksquare\end{array}\right]$};
\node[xscale=0.4, yscale=0.5]  at (1,-2.3) {$\left[\begin{array}{cc|c} \blacksquare & \square & \square \\ \square & \blacksquare & \square \\\hline  \square & \square & \blacksquare \end{array}\right]$};
\node[xscale=0.4, yscale=0.5]  at (3,-2.3) {$\left[\begin{array}{cc|c} \blacksquare & \square & \blacksquare \\ \square & \blacksquare & \square \\ \hline \square & \square & \blacksquare \end{array}\right]$};
\node[xscale=0.4, yscale=0.5]  at (5,-2.3) {$\left[\begin{array}{cc|c} \blacksquare & \square & \square \\ \square & \blacksquare & \blacksquare \\\hline  \square & \square & \blacksquare \end{array}\right]$};
\node[xscale=0.4, yscale=0.5]  at (7,-2.3) {$\left[\begin{array}{cc|c} \blacksquare & \blacksquare &\square\\ \square & \blacksquare &\square \\\hline \square & \square & \blacksquare\end{array}\right]$};
\node[xscale=0.4, yscale=0.5]  at (9,-2.3) {$\left[\begin{array}{cc|c} \blacksquare & \blacksquare &\blacksquare\\ \square & \blacksquare &\square \\\hline \square & \square & \blacksquare\end{array}\right]$};
\node[xscale=0.4, yscale=0.5]  at (11,-2.3) {$\left[\begin{array}{cc|c} \blacksquare & \blacksquare & \blacksquare\\ \square & \blacksquare &\blacksquare \\\hline \square & \square & \blacksquare\end{array}\right]$};
\node[xscale=0.4, yscale=0.5]  at (4.8,-4.5) {$\left[\begin{array}{ccc|c} \blacksquare & \blacksquare & \square& \square\\ \square & \blacksquare &\square & \square\\ \square & \square & \blacksquare& \square \\\hline \square & \square & \square & \blacksquare\end{array}\right]$};
\node[xscale=0.4, yscale=0.5]  at (6.3,-4.5) {$\left[\begin{array}{ccc|c} \blacksquare & \blacksquare &\square& \blacksquare\\ \square & \blacksquare &\square & \square\\\square & \square & \blacksquare& \square\\\hline  \square & \square & \square & \blacksquare\end{array}\right]$};
\node[xscale=0.4, yscale=0.5]  at (7.7,-4.5) {$\left[\begin{array}{ccc|c} \blacksquare & \blacksquare &\square& \blacksquare\\ \square & \blacksquare &\square & \blacksquare\\ \square & \square & \blacksquare& \square\\ \hline\square & \square & \square & \blacksquare\end{array}\right]$};
\node[xscale=0.4, yscale=0.5]  at (9.2,-4.5) {$\left[\begin{array}{ccc|c} \blacksquare & \blacksquare &\square&\square\\ \square & \blacksquare &\square &\square\\ \square & \square & \blacksquare&\blacksquare\\ \hline\square & \square & \square & \blacksquare\end{array}\right]$};
\end{tikzpicture}
\end{center}
\caption{\small{There are $n!$ upper triangular matrix representations of orthogonal $0/1$-simplices.}}
\label{Nfigure13}
\end{figure}  
\begin{rem}{\rm Modulo the action of the hyperoctahedral group, there remain as many as the number of unlabeled trees on $n+1$ vertices. 
Indeed, it is not hard to verify that two matrices $P$ and $R$ representing orthogonal $0/1$-simplices can be transformed into one another 
using operations of type (R),(C) and (X) if and only if the spanning trees of orthogonal edges of the simplices corresponding to $P$ and $R$ are 
isomorphic as graphs.}
\end{rem}

\section{One Neighbor Theorem for a class of nonobtuse simplices}\label{Sect-6}
In this section we will discuss the one neighbor theorem for acute simplices \cite{BrDiHaKr} and generalize it to a larger class of nonobtuse 
simplices. This appears a very nontrivial matter, which can be compared with the complications that arise when generalizing the Perron-
Frobenius theory for positive matrices to nonnegative matrices \cite{BaRa,BePl}.

\subsection{The acute case revisited}
The one neighor theorem for acute $0/1$-simplices reads as follows. We present a alternative proof to the proof in \cite{BrDiHaKr}, based in 
Theorem \ref{th-1}.

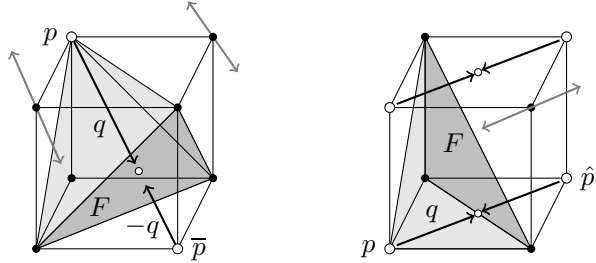
\begin{figure}[h]
\begin{center}
\begin{tikzpicture}[scale=0.94, every node/.style={scale=0.94}]
\draw[fill=gray!20!white] (0,0)--(2.5,1)--(2,2)--(0.5,3)--cycle;
\draw[fill=gray!50!white] (0,0)--(2.5,1)--(2,2)--cycle;
\draw[gray,thick,<->] (-0.35,2.8)--(0.35,1.2);
\draw[gray,thick,<->] (2.15,3.5)--(2.85,2.5);
\node at (2.3,0) {$\overline{p}$};
\node at (0.2,3) {$p$};
\node at (0.9,0.6) {$F$};
\node at (0.9,1.7) {$q$};
\node at (1.5,0.3) {$-q$};
\draw (0,0)--(2,0)--(2,2)--(0,2)--cycle;
\draw (0.5,1)--(2.5,1)--(2.5,3)--(0.5,3)--cycle;
\draw (0,0)--(0.5,1);
\draw (2,0)--(2.5,1);
\draw (0,2)--(0.5,3);
\draw (2,2)--(2.5,3);
\draw (0,0)--(2,2);
\draw (2.5,1)--(0.5,3);
\draw[fill=black] (0,0) circle [radius=0.05];
\draw[fill=black] (2.5,1) circle [radius=0.05];
\draw[fill=black] (2,2) circle [radius=0.05];
\draw[fill=black] (2.5,3) circle [radius=0.05];
\draw[fill=black] (0,2) circle [radius=0.05];
\draw[fill=black] (0.5,1) circle [radius=0.05];
\draw[fill=white] (1.45,1.1) circle [radius=0.05];
\draw[thick,->] (2,0)--(1.55,0.9);
\draw[thick,->] (0.5,3)--(1.4,1.2);
\draw[fill=white] (0.5,3) circle [radius=0.07];
\draw[fill=white] (2,0) circle [radius=0.07];
\begin{scope}[shift={(5,0)}]
\node at (-0.3,0) {$p$};
\node at (2.8,1) {$\hat{p}$};
\draw[fill=gray!20!white] (0,0)--(2,0)--(0.5,3)--cycle;
\draw[fill=gray!50!white] (2,0)--(0.5,3)--(0.5,1)--cycle;
\node at (0.6,0.5) {$q$};
\node at (0.9,1.5) {$F$}; 
\draw (0,0)--(2,0)--(2,2)--(0,2)--cycle;
\draw (0.5,1)--(2.5,1)--(2.5,3)--(0.5,3)--cycle;
\draw (0,0)--(0.5,1);
\draw (2,0)--(2.5,1);
\draw (0,2)--(0.5,3);
\draw (2,2)--(2.5,3);
\draw (2,0)--(0.5,1);
\draw[fill=white] (0,0) circle [radius=0.07];
\draw[fill=black] (0.5,1) circle [radius=0.05];
\draw[fill=black] (2,0) circle [radius=0.05];
\draw[fill=black] (0.5,3) circle [radius=0.05];
\draw[fill=black] (2,2) circle [radius=0.05];
\draw[fill=white] (2.5,1) circle [radius=0.07];
\draw[fill=white] (0,2) circle [radius=0.07];
\draw[fill=white] (2.5,3) circle [radius=0.07];
\draw[fill=white] (1.25,0.5) circle [radius=0.05];
\draw[fill=white] (1.25,2.5) circle [radius=0.05];
\draw[thick,->] (0.1,0.05)--(1.2,0.47);
\draw[thick,->] (2.4,0.95)--(1.3,0.53);
\draw[thick,->] (0.1,2.05)--(1.2,2.47);
\draw[thick,->] (2.4,2.95)--(1.3,2.53);
\draw[gray,thick,<->] (1.3,1.7)--(2.7,2.3);
\end{scope}
\end{tikzpicture}
\end{center}
\caption{\small{Left: for any facet $F$ of an acute $0/1$-simplex, there is only one pair of antipodal vertices $p,\ol{p}$ that may project onto 
$F$. All others end up outside $I^n$ when following the normal $q$ in either direction. Right: in a nonobtuse simplex, there can be more than 
two vertices that remain in $I^n$ when following the normal direction to an interior facet.}}
\label{Nfigure14}
\end{figure}

\begin{Th}[One Neighbor Theorem] \label{th-5} Let $S$ be an acute $0/1$-simplex in $I^n$, and $F$ an $(\nmo)$-facet of $S$ opposite the 
vertex $v$ of $S$. Write $\hat{S}$ for the convex hull of $F$ and $\ol{v}$. Then:\\[3mm]
$\bullet$ $F$ does not lie in an $(\nmo)$-facet of $I^n$;\\[3mm]
$\bullet$ $\hat{S}$ is the only $0/1$-simplex having $F$ as a facet that may be acute, too.\\[3mm]
In words, an acute $0/1$-simplex has at most one acute face-to-face neighbor at each facet.
\end{Th}
{\bf Proof. } Let $q$ be a normal vector to a facet of $F$ opposite a vertex $p$ of an acute $0/1$-simplex $S$. Then due to 
Theorem \ref{th-1}, $q$ has no zero entries. Therefore, the line $v+\alpha q$ parametrized by $\alpha\in\RR$ intersect the interior of $I^n$ if 
and only if $v\in\{p,\ol{p}\}$. Thus, for other vertices $v$ of $I^n$, the altitude from $v$ to the ambient hyperplane of $F$ does not land in 
$I^n$ and thus not in $F$. This is however a necessary condition for the convex hull of $F$ and $v$ to be acute. \hfill $\Box$\\[3mm]
See the left picture in Figure~\ref{Nfigure14} for an illustration of Theorem \ref{th-5} in $I^3$. Only the pair of white vertices end up inside $I^3$ when 
following the normal direction $q$ of the facet $F$. All six other vertices, when projected on the plane containing $F$ end up outside $I^3$, or 
on themselves.

\smallskip

The translation of Theorem \ref{th-5} in terms of linear algebra is as follows. 
\begin{Co} Let $P\in\BB^{n\times(\nmo)}$. The matrix $[P\,|\,v] \in\Bnn$ is a matrix representation of an acute $0/1$-simplex for at most 
one pair of antipodal points $v\in\{p,\ol{p}\}\subset\BB^n$.
\end{Co}
The one neighbor theorem dramatically restricts the number of $0/1$-polytopes that can be face-to face triangulated by acute simplices. For 
instance, only from dimension $n=7$ onwards there exists a pair of face-to-face acute simplices in $I^n$. In $I^7$ it is the Hadamard regular 
simplex \cite{Gr} and its face-to-face neighbor, which is unique modulo the action of the hyperoctahedral group, and which has the one-but-
largest volume in $I^7$ over all acute $0/1$-simplices \cite{BrCi2}. Also in \cite{BrDiHaKr} the theorem turned out useful in constructing all 
possible face-to-face triangulations of $I^n$ consisting on nonobtuse simplices only, due to the following sharpening of the statement.

\begin{Co} Each acute $0/1$-simplex $S$ in $I^n$ has at most one face-to-face nonobtuse neighbor at each of its facets.
\end{Co}
{\bf Proof.} This follows from the fact that in the proof of Theorem \ref{th-5}, the altitudes from $v\not=\{p,\ol{p}\}$ intersect $I^n$ only in 
$v$ itself.\hfill $\Box$\\[3mm]
A natural question is what can be proved for nonobtuse-$0/1$ simplices. Theorem \ref{th-2} showed that a normal to a facet of a nonobtuse 
$0/1$-simplex $S$ can have entries equal to zero. Writing $\zeros(q)$ for the number of entries of $q$ equal to zero, there 
are$2^{\zeros(q)+1}$ vertices $v$ of $I^n$ from which the altitudes starting at $v$ onto the plane containing $F$ do not leave $I^n$. This is 
illustrated in the right picture in Figure~\ref{Nfigure14}. The normal vector $q$ to the facet $F$ has one zero entry: $\zeros(q)=1$. The altitudes from the 
$2^2$ white vertices of $I^3$ onto the plane containing $F$ lie in $I^3$.\\[3mm] 
Nevertheless, only the altitudes from $p$ and $\hat{p}$ land on $F$ itself, and we see that the interior facet $F$ of $S$ in $I^3$ has exactly 
one nonobtuse neighbor. It is tempting to conjecture that the one neighbor theorem holds also for nonobtuse $0/1$-simplices. The only 
adaption to make is then, based on the example in Figure~\ref{Nfigure14}, that instead of the antipodal $\ol{p}$ of $p$ in $I^n$, the second vertex 
$\hat{p}$ such that the convex hull of $F$ with $\hat{p}$ is a nonobtuse simplex would satisfy 
\be\label{restr-ant} \hat{p}_j=1-p_j \hdrie \Leftrightarrow \hdrie q_j\not=0 \und \hat{p}_j = 0 \hdrie \Leftrightarrow \hdrie q_j=0. \ee
In words, $\hat{p}$ is the antipodal of $p$ restricted to the $(n-\zeros(q))$-dimensional $n$-cube facet that contains both $p$ and $q$. In 
Figure~\ref{Nfigure14}, $\hat{p}$ is the antipodal of $p$ in the bottom square facet of $I^3$.\\[3mm]
As an attempt to prove this conjecture, one may try to demonstrate that the remaining white vertices in the top square facet of $I^3$, although 
their altitudes lie in $I^3$, cannot fall onto $F$. Although we did not succeed in doing so, the nonobtusity of $S$ is a necessary condition. To 
see this, consider the $0/1$-simplex $S$ in $I^5$ with matrix representation
\be \left[\begin{array}{ccccc} 1 & 1 & 1 & 0 & 0\\ 1 & 1 & 0 & 1 & 0\\ 1 & 0 & 0 & 0 & 1\\ 0 & 0 & 1 & 1 & 0 \\ 0 & 1 & 1 & 0 & 1\end{array}\right] \hdrie \hdrie\mbox{\rm with} \hdrie q = \frac{1}{2}\left[\begin{array}{c} 0 \\ 1 \\ 1 \\ 1 \\ 1\end{array}\right],\ee
and $q$ is the normal to the facet $F$ of $S$ opposite the origin. Observe that the line from the origin to the vertex $2q$ in $I^5$ intersects 
$F$ in the midpoint of its edge between the two vertices of $S$ in the last two columns of $P$. This shows that both the origin and its 
antipodal in the bottom $4$-facet of $I^5$ as defined in (\ref{restr-ant}) land in $F$ when following their respective altitudes. But also the line 
between $e_1^5$ and $e^5$ intersects $F$ in the midpoint of its edge between the two vertices in the first and third column of $P$. Thus, 
also both the vertices $e_1^5$ and $e^5$ in the top $4$-facet of $I^5$ land in $F$ when following their altitudes. Thus, for a $0/1$-simplex 
that is not nonobtuse, it can occur that more than two vertices of $I^n$ project orthogonally onto $F$.

\subsection{More on fully indecomposable nonobtuse simplices}
In Section~\ref{Sect-6.3} we will study the one neighbor theorem in the context of partly decomposable 
nonobtuse simplices. For this, but also for its own interest, we derive two auxiliary results on fully indecomposable nonobtuse simplices.

\begin{Le}\label{lem-6} Each representation $P\in\Bnn$ of a fully indecomposable $0/1$-simplex $S$ satisfies
\be P^\top P  \geq I+e^n\left(e^n\right)^\top. \ee
In geometric terms this implies that all triangular facets of $S$ are acute.
\end{Le}
{\bf Proof. } A standard type of argument is the following. Write $D$ for the diagonal matrix having the same diagonal 
entries as $B=(P^\top P)^{-1}$, and let $C=D-B$. Then $C\geq 0$, and
\be B = D-C = D(I-D^{-1}C) \und P^\top P =B^{-1} =  (I-D^{-1}C)^{-1}D^{-1}. \ee
Because $B$ is an M-matrix \cite{Joh,JoSm}, the spectral radius of $D^{-1}C$ is less than one, and the following Neumann series converges:
\be (I-D^{-1}C)^{-1} = \sum_{j=0}^\infty \left(D^{-1}C\right)^j. \ee
Since $P$ is fully indecomposable, $P^\top P$ is irreducible, and thus $B$ is irreducible. But then, so are $C$ and $D^{-1}C$. 
Because of the latter, for each pair $k,\ell$ there is a $j$ such that $e_k^\top(D^{-1}C)^j e_\ell >0$. This proves that 
$B^{-1}=P^\top P>0$. Since its entries are integers, $P^\top P \geq e^n(e^n)^\top$. Thus, each pair of edges of $S$ that meet 
at the origin makes an acute angle. As by Theorem \ref{th-1} all matrix representations of $S$ are fully indecomposable, we conclude 
that all triangular facets of $S$ are acute. This implies that any diagonal entry of $P^\top P$ is greater than the remaining entries in 
the same row. Indeed, if two entries in the same row would be equal, then $p_j^\top (p_j-p_i)=0$, which corresponds to two edges of $S$ making a right angle.\hfill $\Box$

\begin{Co}\label{co-5} Let $P$ be a fully indecomposable matrix representation of a nonobtuse $0/1$-simplex. 
If $\hat{P}$ equals $P$ with one column replaced by its antipodal, then $\hat{P}^\top\hat{P}>0$.
\end{Co}
{\bf Proof. } Without loss of generality, assume that
\be P = [p\,|\,P_1] \und \hat{P} = [\ol{p}\,|\,P_1] \ee
with $p\in\BB^n$. Then
\be \hat{P}^\top\hat{P} = \left[\begin{array}{cc} \ol{p}^\top\ol{p} & \ol{p}^\top P_1 \\ P_1\ol{p} & P_1^\top P_1 \end{array}\right], \ee
and $P_1^\top P_1>0$ because $P^\top P>0$ as proved in Lemma \ref{lem-5}. Due to the fact that for all $a,b\in\BB^n$,
\be a^\top \ol{b} = a^\top (e^n-b) = a^\top(a-b),\ee
we see that also $P_1^\top\ol{p}>0$. Indeed, a zero entry would contradict that the diagonal entries of $P^\top P$ are greater than its 
off-diagonal entries, as proved in Lemma \ref{lem-6}.\hfill $\Box$
  
\begin{rem}{\rm The matrix $\hat{P}$ in Corollary \ref{co-5} is not always fully indecomposable. See
\be P=\left[\begin{array}{ccc} 0 & 1 & 1 \\ 1 & 0 & 1\\ 1 & 1 & 0\end{array}\right] \und \hat{P}=\left[\begin{array}{ccc} 1 & 1 & 1 \\ 0 & 0 & 1\\ 0 & 1 & 0\end{array}\right], \ee 
where the first columns of both matrices are each other's antipodal. This example also shows that the top left diagonal entry of 
$\hat{P}^\top\hat{P}$ need not be greater than the other entries in its row.}
\end{rem}
We end this section with a theorem which was proved by inspection of a finite number of cases. We refer to \cite{BrCi2} for details 
on how to computationally generate the necessary data.

\begin{Th}\label{th-6} Each fully indecomposable nonobtuse $0/1$-simplex in  $\SS^n$ with $n\leq 8$ is acute. There exist fully 
indecomposable nonobtuse $0/1$-simplices in $\SS^n$ with $n\geq 9$ that are not acute.
\end{Th}
{\bf Proof. } See \cite{BrCi2} for details on an algorithm to compute $0/1$-matrix representations of $0/1$-simplices modulo the action 
of the hyperoctahedral group. By inspection of all $0/1$-simplices of dimensions less than or equal to $8$, we conclude the first 
statement. For the second statement, we give an example. The $9\times 9$ matrix in Figure~\ref{Nfigure15} represents a nonobtuse simplex 
that is not acute, but $P$ is fully indecomposable. \hfill $\Box$
\begin{figure}[h]
\begin{center}
\begin{tikzpicture}[scale=0.94, every node/.style={scale=0.94}]
\node[scale=0.7] at (0,0)  {$P = \left[\begin{array}{rrrrrrrrr}
  1  & 1 &  0 &  0 &  1 &  1 &  1 &  1 &  0\\
  1  & 0 &  1 &  1 &  1 &  0 &  0 &  1 &  1\\ 
  1  & 0 &  1 &  1 &  0 &  1 &  1 &  0 &  1\\
  0  & 1 &  1 &  1 &  1 &  0 &  1 &  0 &  1\\
  0  & 1 &  1 &  1 &  0 &  1 &  0 &  1 &  1\\
  0  & 0 &  1 &  1 &  1 &  1 &  1 &  1 &  0\\
  0  & 0 &  1 &  0 &  1 &  1 &  0 &  0 &  1\\
  0  & 0 &  1 &  0 &  0 &  0 &  1 &  1 &  1\\
  0  & 0 &  0 &  1 &  1 &  1 &  1 &  1 &  1
  \end{array}\right], \hspace{3mm}P^{-\top}= \dfrac{1}{20} 
   \left[\begin{array}{rrrrrrrrr}   
     6  &   6  &  -2  &  -6  &   2  &   2  &   2  &   2  &  -2\\
     7  &  -3  &   1  &   3  &   4  &  -6  &  -6  &   4  &   1\\
     7  &  -3  &   1  &   3  &  -6  &   4  &   4  &  -6  &   1\\
    -3  &   7  &   1  &   3  &   4  &  -6  &   4  &  -6  &   1\\
    -3  &   7  &   1  &   3  &  -6  &   4  &  -6  &   4  &   1\\
    -4  &  -4  &   8  &   4  &   2  &   2  &   2  &   2  & -12\\
    -2  &  -2  &   4  &  -8  &   6  &   6  &  -4  &  -4  &   4\\
    -2  &  -2  &   4  &  -8  &  -4  &  -4  &   6  &   6  &   4\\
    -4  &  -4  & -12  &   4  &   2  &   2  &   2  &   2  &   8
    \end{array}\right], \hspace{3mm}\dfrac{1}{20} 
\left[\begin{array}{r} 10 \\ 5 \\ 5 \\ 5 \\ 5 \\ 0 \\ 0 \\ 0 \\ 0 \end{array}\right]$};
\end{tikzpicture}
\end{center}
\caption{\small{Example of a fully indecomposable matrix representation $P$ of a nonobtuse simplex $S$ that is not acute. The vector 
on the right is its normal $q$ to the facet $F$ opposite the origin. Since $q$ has entries equal to zero, $S$ cannot be acute. But 
$(P^\top P)^{-1}$ satisfies (\ref{eq-1.1}) and (\ref{eq-1.2}), hence $S$ is nonobtuse. Note that none of the other normals has a zero entry.}}
\label{Nfigure15}
\end{figure}  

\smallskip

Thus, Theorem \ref{th-6} and Figure~\ref{Nfigure14} prove that the {\em indecomposability} of a matrix representation of a nonobtuse $0/1$-simplex $S$ is, in 
fact, a {\em weaker} property than the {\em acuteness} of $S$. \\[3mm]
Especially since the two concepts coincide up to dimension eight, this came as a surprise. Citing G\"unther Ziegler in 
Chapter 1 of {\em Lectures on $0/1$-Polytopes} \cite{KaZi}: ``{\em Low-dimensional intuition does not work!}\,''. 
See \cite{BrCi2} for more such examples in the context of $0/1$-simplices.

\subsection{A One Neighbor Theorem for partly decomposable simplices}\label{Sect-6.3}
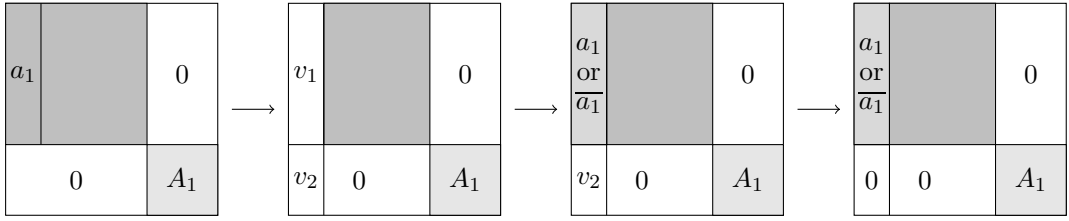
\begin{figure}[h]
\begin{center}
\begin{tikzpicture}[scale=0.94, every node/.style={scale=0.94}]
\draw (0,0)--(3,0)--(3,3)--(0,3)--cycle;

\draw[fill=gray!20!white] (2,0)--(3,0)--(3,1)--(2,1)--cycle;
\draw[fill=gray!50!white] (0,1)--(2,1)--(2,3)--(0,3)--cycle;

\draw (0.5,1)--(0.5,3);

\node at (2.5,0.5) {$A_1$};
\node at (1,0.5) {$0$};
\node at (2.5,2) {$0$};
\node at (0.25,2) {$a_1$};

\draw[->] (3.2,1.5)--(3.8,1.5);

\begin{scope}[shift={(4,0))}]

\draw (0,0)--(3,0)--(3,3)--(0,3)--cycle;

\draw[fill=gray!20!white] (2,0)--(3,0)--(3,1)--(2,1)--cycle;
\draw[fill=gray!50!white] (0.5,1)--(2,1)--(2,3)--(0.5,3)--cycle;

\draw (0.5,0)--(0.5,3);
\draw (0,1)--(0.5,1);

\node at (2.5,0.5) {$A_1$};
\node at (1,0.5) {$0$};
\node at (2.5,2) {$0$};

\node at (0.25,2) {$v_1$};
\node at (0.25,0.5) {$v_2$};

\draw[->] (3.2,1.5)--(3.8,1.5);

\end{scope}

\begin{scope}[shift={(8,0))}]

\draw (0,0)--(3,0)--(3,3)--(0,3)--cycle;

\draw[fill=gray!20!white] (2,0)--(3,0)--(3,1)--(2,1)--cycle;
\draw[fill=gray!50!white] (0.5,1)--(2,1)--(2,3)--(0.5,3)--cycle;
\draw[fill=gray!30!white] (0,1)--(0.5,1)--(0.5,3)--(0,3)--cycle;

\draw (0.5,0)--(0.5,3);
\draw (0,1)--(0.5,1);

\node at (2.5,0.5) {$A_1$};
\node at (1,0.5) {$0$};
\node at (2.5,2) {$0$};

\node at (0.25,2.4) {$a_1$};
\node at (0.25,2) {or};
\node at (0.25,1.6) {$\overline{a_1}$};
\node at (0.25,0.5) {$v_2$};

\draw[->] (3.2,1.5)--(3.8,1.5);

\end{scope}

\begin{scope}[shift={(12,0))}]

\draw (0,0)--(3,0)--(3,3)--(0,3)--cycle;

\draw[fill=gray!20!white] (2,0)--(3,0)--(3,1)--(2,1)--cycle;
\draw[fill=gray!50!white] (0.5,1)--(2,1)--(2,3)--(0.5,3)--cycle;
\draw[fill=gray!30!white] (0,1)--(0.5,1)--(0.5,3)--(0,3)--cycle;

\draw (0.5,0)--(0.5,3);
\draw (0,1)--(0.5,1);

\node at (2.5,0.5) {$A_1$};
\node at (1,0.5) {$0$};
\node at (2.5,2) {$0$};

\node at (0.25,2.4) {$a_1$};
\node at (0.25,2) {or};
\node at (0.25,1.6) {$\overline{a_1}$};
\node at (0.25,0.5) {$0$};

\end{scope}
\end{tikzpicture}
\end{center}
\caption{\small{Steps in proving a one neighbor theorem for partly decomposable nonobtuse simplices with two fully indecomposable 
blocks representing acute simplices.}}
\label{Nfigure16}
\end{figure}

Let $S$ be a partly decomposable nonobtuse simplex. Then according to Corollary \ref{co-3}, $S$ has a matrix representation
\be P = \left[\begin{array}{cc} N & 0 \\ 0 & A\end{array}\right] \ee 
in which $A$ is fully indecomposable. We will discuss some cases in which a modified version of the One Neighbor Theorem \ref{th-5} 
holds also for nonobtuse simplices that are not acute.\\[3mm]
{\bf Case I.} To illustrate the main line of argumentation, consider first the simplest case, which is that $S\in\SS^n$ can be represented by
\be\label{eq-21} P = \left[\begin{array}{cc} A_2 & 0 \\ 0 & A_1\end{array}\right],\hdrie\mbox{\rm with} \hdrie A_1\in\BB^{k\times k} \und A_2\in\BB^{(n-k)\times(n-k)} \ee
and in which $A_1$ and $A_2$ represent acute simplices $S_1$ and $S_2$. Then $A_1$ and $A_2$ are fully indecomposable, 
and $S_1$ and $S_2$ satisfy the One Neighbor Theorem \ref{th-5}. Assume first that neither $A_1$ or $A_2$ equals the $1\times 1$ matrix $[\,1\,]$.\\[3mm] 
{\bf Notation. } We will write $X^j(y)$ for the matrix $X$ with column $j$ replaced by $y$.\\[3mm]
Now, let $v\in\BB^n$, partitioned as $v^\top = (v_1^\top \,\,v_2^\top)$ with $v_1\in\BB^k$. Assume that the block lower 
triangular matrix $P^1(v)$ represents a nonobtuse simplex. This implies that its top-left diagonal block $A_2^1(v_1)$ does 
so, too, hence $v_1=a_1=Ae_1^k$ or $v_1=\ol{a_1}$ by Theorem \ref{th-5}. Corollary \ref{co-5} shows that in both cases 
$A_2^1(v_1)^\top A_2^1(v_1) > 0$. But then Theorem \ref{th-3} in combination with the observations in Remarks \ref{rem-4} 
and \ref{rem-6} proves that $v_2=0$, because all columns in the off-diagonal block must be copies of one and the same column 
of $A_1$. Hence, there is at most one $v\in\BB^n$ other than $Pe_1^n$ such that $P^1(v)$ is nonobtuse. See Figure~\ref{Nfigure16} for a sketch of the proof.

\smallskip

Clearly, the same argument can be applied to prove that for all $j\in\{1,\dots,n\}$, the matrix $P^j(v)$ represents a nonobtuse $0/1$-simplex 
for at most one $v\in\BB^n$ other than $Pe_j^n$. This proves that for each column $p$ of $P$, the facet $F_p$ of $S$ opposite $p$ has at 
most one nonobtuse neighbor. It remains to prove the same for the facet $F_0$ of $S$ opposite the origin. But because $S_1$ and $S_2$ are 
by assumption acute, the normals $q_1$ and $q_2$ to their respective facets opposite the origin, which satisfy $A_1^\top q_1 = e_k^k$ and 
$A_2^\top q_2 = e_{n-k}^{n-k}$, are both positive. But then so is the normal $q$ of $F_0$, which satisfies $P^\top q = e_n^n$, and 
hence $q^\top = (q_1^\top\,\,q_2^\top)>0$. And thus, apart from the origin, only $e_n^n$ can form a nonobtuse simplex together with $F_0$.

\begin{rem}{\rm Note that this last argument does not hold if $A_1$ and $A_2$ are merely assumed to represent fully indecomposable 
nonobtuse simplices: the $9\times 9$ matrix in Figure~\ref{Nfigure15} shows that the normal of the facet opposite the origin may contain entries equal to zero.}
\end{rem}
To finish the case in which $S$ has a matrix representation as in (\ref{eq-21}), assume without loss of generality that $A_1=[\,1\,]$ and 
$A_2\not=[\,1\,]$.  Then the facet $F$ of $S$ opposite the last column of $P$ lies in a cube facet, and thus it cannot have a nonobtuse 
face-to-face neighbor. For the remaining $n$ facets of $S$, arguments as above apply, and we conclude that $S$ has at most one nonobtuse 
neighbor at each of its facets. The remaining case that $A_1=A_2=[\,1\,]$ is trivial.\\[3mm]
Note that if a nonobtuse $0/1$-simplex has a block diagonal matrix representation with $p>2$ blocks, each representing an acute simplex, the 
result remains valid, based on a similar proof.\\[3mm]
{\bf Case II.} Assume now that the matrix representation $P$ of a nonobtuse $0/1$-simplex $S$ has the form
\be\label{case-iii} P = \left[\begin{array}{cc} N_1 & 0 \\ 0 & A_1\end{array}\right], \ee 
where $A_1\in\BB^{(n-k)\times(n-k)}$ represents an acute simplex and $N_1$ a merely nonobtuse simplex. Using similar 
arguments as in Case I it is easily seen that the only two choices of $v\in\BB^n$ such that $P^j(v)$ with $k+1\leq j \leq n$ is nonobtuse are
\be v = Pe_j^n = \left[\begin{array}{c} 0 \\ a_j \end{array}\right] \und v = \left[\begin{array}{c} 0 \\ \ol{a_j} \end{array}\right], \ee 
as no additional properties of $N_1$ need to be known. This changes if we examine the matrix $P^j(v)$ with $1\leq j\leq k$, as it is generally not 
true that $N_1^\top N_1>0$. A way out is the following. Assume that also $N_1$ is partly decomposable, then using only row and column permutations, 
we can first transform $N_1$ into the form
\be N_1 \overset{(C)+(R)}{\longrightarrow } \left[\begin{array}{cc} N_2 & R \\ 0  & A_2\end{array}\right] \ee 
where we assume that $A_2$ represents an acute simplex. Then reflecting the vertex to the origin such that the block above $A_2$ becomes zero, we find that
\be\label{eq-22} P \sim  \tilde{P} = \left[\begin{array}{c|c|c}  N_2 & R & 0 \\ \hline 0  & A_2 & 0 \\\hline 0 & 0 & A_1\end{array}\right] \sim  \left[\begin{array}{c|c|c}  \tilde{N}_2 & 0 & \tilde{R} \\\hline 0  & A_2 & 0 \\\hline 0 & 0 & A_1\end{array}\right] = \hat{P},\ee
where $\tilde {R}$ has the same columns as $R$ but possibly a different number of them. Now, select a column of $\hat{P}$ that contains entries of $A_2$ and 
replace it by $v$, partitioned as $v^\top = (v_1^\top \,\,v_2^\top\,\,v_3^\top)$. Because the bottom right $2\times 2$ block part of $\hat{P}$ is a matrix 
representation of a nonobtuse simplex as considered in Case I, we conclude that $v_3=0$. Because the top left $2\times 2$ block part of $\tilde{P}$ is a 
matrix representation as in (\ref{case-iii}), we conclude that $v_3=0$ and $v_2$ is a column of $A_2$ or its antipodal. Thus, also the facets of the 
vertices of $S$ corresponding to its indecomposable part $A_2$ all have at most one nonobtuse neighbor.\\[3mm] 
Now, this process can be inductively repeated in case $\tilde{N_2}$ is partly decomposable with a fully indecomposable part that represents an acute 
simplex, and so on, until a fully indecomposable top left block $A_p$ remains. This block presents vertices for which it still needs to be proved that 
their opposite facets have at most one nonobtuse neighbor. To illustrate how to do this, consider the case $p=3$. Or, in other words, assume that 
$N_2$ in (\ref{eq-22}) represents an acute simplex. Replace one of the corresponding columns of $\tilde{P}$ by $v$ partitioned as 
$v^\top = (v_1^\top \,\,v_2^\top\,\,v_3^\top)$. Then $v_3=0$ because the $(1,3)$ block of $\tilde{P}$ equals zero. Similarly, because 
the $(1,2)$-block in $\hat{P}$ equals zero, we find that $v_2=0$. And thus, $v_3$ is a column of $N_2$ or its antipodal. For $p>3$, we 
can do the same one by one for the blocks at positions $(1,p),\dots,(1,2)$.\\[3mm]
The analysis in this section can be summarized in the following theorem. 

\begin{Th} Let $S$ be a nonobtuse $0/1$-simplex whose fully indecomposable components are all acute. Then $S$ has at most one face-to-face neighbor at each  of its interior facets.
\end{Th}

\subsection*{Acknowledgments}
Jan Brandts and Apo Cihangir acknowledge the support by Research Project 613.001.019 of the Netherlands Organisation for Scientific Research (NWO), and are grateful to Michal K\v{r}\'{\i}\v{z}ek for comments and discussions on earlier versions of the manuscript.

\end{document}